\documentclass[11pt]{article}
\usepackage{amsmath}
\usepackage{amssymb}

\usepackage{theorem}
\numberwithin{equation}{section}

\newtheorem{theorem}{Theorem}[section]
\newtheorem{proposition}[theorem]{Proposition}

\newtheorem{lemma}[theorem]{Lemma}
\newtheorem{definition}{Definition}[section]

\begin{document}
\title{Almost Morawetz estimates and global well-posedness for the defocusing $L^2$-critical nonlinear Schr{\"o}dinger equation in higher dimensions}
\date{\today}
\author{Ben Dodson}
\maketitle

\noindent \textbf{Abstract:} In this paper, we consider the global well-posedness of the defocusing, $L^{2}$ - critical nonlinear Schr{\"o}dinger equation in dimensions $n \geq 3$. Using the I-method, we show the problem is globally well-posed in $n = 3$ when $s > \frac{2}{5}$, and when $n \geq 4$, for $s > \frac{n - 2}{n}$. We combine energy increments for the I-method, interaction Morawetz estimates, and almost Morawetz estimates to prove the result.

\section{Introduction}
The defocusing, $L^{2}$ - critical nonlinear Schr{\"o}dinger equation

\begin{equation}\label{0.1}
\aligned
i u_{t} + \Delta u &= |u|^{4/n} u, \\
u(0,x) &= u_{0}(x) \in H^{s}(\mathbf{R}^{n}),
\endaligned
\end{equation}

\noindent has a local solution on some interval $[0, T]$, $T(\| u_{0} \|_{H^{s}(\mathbf{R}^{n})}) > 0$ when $s > 0$. (See \cite{CaWe1}.) $(\ref{0.1})$ also has a local solution when $u_{0} \in L^{2}(\mathbf{R}^{n})$  on $[0, T)$, $T(u_{0}) > 0$, where $T$ depends on the profile of the initial data, not just its size. For global well-posedness to fail, and a solution to $(\ref{1.1})$ only exist on a maximal interval $[0, T_{\ast})$, $T_{\ast} < \infty$, then

\begin{equation}\label{0.1.1}
\lim_{t \rightarrow T_{\ast}} \| u(t) \|_{H^{s}(\mathbf{R}^{n})} = \infty
\end{equation}

\noindent for all $s > 0$. $(\ref{0.1})$ has the conserved quantities:

\begin{equation}\label{0.2}
M(u(t)) = \int |u(t,x)|^{2} dx = M(u(0)),
\end{equation}

\begin{equation}\label{0.3}
E(u(t)) = \frac{1}{2} \int |\nabla u(t,x)|^{2} dx + \frac{n}{2n + 4} \int |u(t,x)|^{2 + 4/n} dx.
\end{equation}

\noindent Combining the fact that $E(u(t))$ is positive definite, the Sobolev embedding theorem, and $(\ref{0.1.1})$; \cite{CaWe} proved $(\ref{0.1})$ is globally well-posed when $u_{0} \in H^{1}(\mathbf{R}^{n})$.\vspace{5mm}

\noindent Furthermore, a solution to the equation

\begin{equation}\label{0.3.1}
\aligned
i u_{t} + \Delta u &= |u|^{\alpha} u, \\
u(0,x) &= u_{0}(x),
\endaligned
\end{equation}

\noindent can be rescaled in the following manner. If $u(t,x)$ is a solution to $(\ref{0.3.1})$ on $[0, T_{0}]$, then $$u_{\lambda}(t,x) = \frac{1}{\lambda^{2/\alpha}} u(\frac{t}{\lambda^{2}}, \frac{x}{\lambda})$$ is a solution to $(\ref{0.3.1})$ on $[0, \lambda^{2} T_{0}]$ with initial data $$\frac{1}{\lambda^{2/\alpha}} u_{0}(\frac{x}{\lambda}).$$ $\alpha = \frac{4}{n}$ is called the $L^{2}$-critical exponent because a brief calculation will show that when $\alpha = \frac{4}{n}$,

\begin{equation}\label{0.3.2}
\| u_{\lambda} \|_{L_{t}^{\infty} L_{x}^{2}([0, \lambda^{2} T_{0}] \times \mathbf{R}^{n})} = \| u \|_{L_{t}^{\infty} L_{x}^{2}([0, T_{0}] \times \mathbf{R}^{n})}.
\end{equation}

\noindent Indeed, for any $n$-admissible pair $(p,q)$,

\begin{equation}\label{0.3.3}
\| u_{\lambda} \|_{L_{t}^{p} L_{x}^{q}([0, \lambda^{2} T_{0}] \times \mathbf{R}^{n})} = \| u \|_{L_{t}^{p} L_{x}^{q}([0, T_{0}] \times \mathbf{R}^{n})}.
\end{equation}

\noindent (Admissible pairs will be discussed in greater detail in $\S 2$.)\vspace{5mm}

\noindent Many have endeavored to prove global well-posedness for less regular data, $u_{0} \in H^{s}(\mathbf{R}^{n})$, $s < 1$. The first progress was made in \cite{B}, proving global well-posedness for $s > \frac{3}{5}$ when $n = 2$ via the Fourier truncation method. In addition to proving global well-posedness, \cite{B} proved

\begin{equation}\label{0.4}
u(t,x) - e^{it \Delta} u_{0} \in H^{1}(\mathbf{R}^{2}).
\end{equation}

\noindent This method was modified in \cite{CKSTT1} to produce the I - method, proving global well-posedness of $(\ref{0.1})$ when $n = 2$, $s > \frac{4}{7}$. (\cite{CKSTT1} also discussed the cubic nonlinear Schr{\"o}dinger equation when $n = 3$, but that equation will not be discussed here, as it is $\dot{H}^{1/2}$ - critical.)\vspace{5mm}

\noindent Since then, several improvements have been made when $n = 2$. In particular, improvements have utilized an almost Morawetz estimate. (See \cite{CGT}, \cite{CR}.) Currently, the best known result for $n = 2$ is

\begin{theorem}\label{t0.1}
$(\ref{0.1})$ is globally well-posed when $n = 2$ for $s > \frac{1}{4}$.
\end{theorem}

\noindent \emph{Proof:} See \cite{D}.\vspace{5mm}

\noindent In \cite{DPST}, the I-method was extended to prove global well-posedness results for $(\ref{0.1})$ when $n \geq 3$. The chief difficulty with extending to $n \geq 3$ is that the nonlinearity $|u|^{4/n} u$ is no longer "algebraic" when $n > 2$. That is, $|u|^{4/n} u$ is no longer a polynomial of $u$ and $\bar{u}$ when $n > 2$. Nevertheless, it was proved that

\begin{theorem}\label{t0.2}
$(\ref{0.1})$ is globally well-posed for $s > \frac{\sqrt{7} - 1}{3}$ when $n = 3$, and $s > \frac{-(n - 2) + \sqrt{(n - 2)^{2} + 8(n - 2)}}{4}$ when $n \geq 4$.
\end{theorem}

\noindent \emph{Proof:} See \cite{DPST}.\vspace{5mm}

\noindent In this paper, we will prove

\begin{theorem}\label{t0.3}
When $n \geq 4$, $(\ref{0.1})$ is globally well-posed for $u_{0} \in H^{s}(\mathbf{R}^{n})$, $s > \frac{n - 2}{n}$. Moreover,

\begin{equation}\label{0.4.1}
\sup_{t \in [0, T_{0}]} \| u(t) \|_{H^{s}(\mathbf{R}^{n})} \leq C(\| u_{0} \|_{H^{s}(\mathbf{R}^{n})}) T_{0}^{\frac{(n - 2)(1 - s)^{2}}{2(ns - (n - 2))}}.
\end{equation}
\end{theorem}

\noindent \vspace{5mm}

\begin{theorem}\label{t0.4}
When $n = 3$, $(\ref{0.1})$ is globally well-posed for $u_{0} \in H^{s}(\mathbf{R}^{n})$, $s > \frac{n - 2}{n}$. Moreover,

\begin{equation}\label{0.4.2}
\sup_{t \in [0, T_{0}]} \| u(t) \|_{{H}^{s}(\mathbf{R}^{3})} \leq C(\| u_{0} \|_{H^{s}(\mathbf{R}^{3})}) T_{0}^{\frac{(1 - s)}{5s - 2}+}.
\end{equation}
\end{theorem}\vspace{5mm}

\noindent \textbf{Description of Method:}

\noindent For $u_{0} \in H^{s}(\mathbf{R}^{n})$, $s < 1$, the I - operator is defined to be the Fourier multiplier

\begin{equation}\label{0.5}
m(\xi) = \left\{
           \begin{array}{ll}
             1, & \hbox{$|\xi| \leq N$;} \\
             (\frac{N}{|\xi|})^{1 - s}, & \hbox{$|\xi| > N$.}
           \end{array}
         \right.
\end{equation}

\noindent Then if $u(t,x)$ solves $(\ref{0.1})$, $Iu(t,x)$ solves

\begin{equation}\label{0.6}
\aligned
i Iu_{t} + I \Delta u = I(|u|^{4/n} u), \\
Iu(0,x) \in H^{1}(\mathbf{R}^{n}).
\endaligned
\end{equation}

\begin{equation}\label{0.6.1}
\aligned
\| Iu \|_{H^{1}(\mathbf{R}^{n})} &\lesssim N^{1 - s} \| u \|_{H^{s}(\mathbf{R}^{n})}, \\
\| u \|_{H^{s}(\mathbf{R}^{n})} &\lesssim \| Iu \|_{H^{1}(\mathbf{R}^{n})},
\endaligned
\end{equation}

\noindent so $E(Iu(t))$ very effectively controls $\| u \|_{H^{s}(\mathbf{R}^{n})}$. The chief difficulty is that, unlike $(\ref{0.3})$, $E(Iu(t))$ is not a conserved quantity, rather,

\begin{equation}\label{0.7}
\frac{d}{dt} E(Iu(t)) = Re \int (\overline{Iu_{t}(t,x)}) (I(|u(t,x)|^{4/n} u(t,x)) - |Iu(t,x)|^{4/n} Iu(t,x)) dx
\end{equation}

\begin{equation}\label{0.8}
= -Im \int \overline{I \Delta u(t,x)} (I(|u(t,x)|^{4/n} u(t,x)) - |Iu(t,x)|^{4/n} Iu(t,x)) dx
\end{equation}

\begin{equation}\label{0.9}
+ Im \int \overline{I(|u|^{4/n} u)} (I(|u(t,x)|^{4/n} u(t,x)) - |Iu(t,x)|^{4/n} Iu(t,x)) dx.
\end{equation}

\noindent To get around the fact that $|u|^{4/n} u$ is not algebraic, we use the fact that $|Iu|^{4/n}$ can be very effectively approximated by $I(|u|^{4/n})$. Therefore, the analysis of $(\ref{0.8})$ and $(\ref{0.9})$ can be split into the analysis of a "main term" and also a "remainder term", yielding a smaller energy increment than in \cite{DPST}. For the purposes of this paper, $(\ref{0.8})$ will be called the linear term, and $(\ref{0.9})$ will be called the nonlinear term.\vspace{5mm}

\noindent In $\S 2$, some preliminary results from harmonic analysis will be discussed. In $\S 3$ some of the smoothness properties of $|u|^{4/n}$ and $|Iu|^{4/n}$ that will be needed later will be proved. In $\S 4$, we will prove

\begin{equation}\label{0.10}
|\int_{t_{1}}^{t_{2}} (\ref{0.8}) dt| \lesssim \frac{1}{N^{4/n-}} \| \langle \nabla \rangle Iu \|_{S^{0}([t_{1}, t_{2}] \times \mathbf{R}^{n})}^{2 + 4/n},
\end{equation}

\noindent when $n \geq 4$. In $\S 5$, using a slightly different method, we will prove

\begin{equation}\label{0.11}
|\int_{t_{1}}^{t_{2}} (\ref{0.8}) dt| \lesssim \frac{1}{N^{1-}} \| \langle \nabla \rangle Iu \|_{S^{0}([t_{1}, t_{2}] \times \mathbf{R}^{n})}^{10/3},
\end{equation}

\noindent when $n = 3$. In $\S 6$ we will prove

\begin{equation}\label{0.12}
|\int_{t_{1}}^{t_{2}} (\ref{0.9}) dt| \lesssim \| \langle \nabla \rangle Iu \|_{S^{0}([t_{1}, t_{2}] \times \mathbf{R}^{n})}^{2 + 4/n} \left\{
                             \begin{array}{ll}
                               \frac{1}{N^{1-}}, & \hbox{when $n = 3$;} \\
                               \frac{1}{N^{\frac{4}{n}-}}, & \hbox{when $n \geq 4$.}
                             \end{array}
                           \right.
\end{equation}

\noindent In $\S 7$ we will prove Theorem $\ref{t0.3}$. In $\S 8$ we will prove an almost Morawetz estimate for $n = 3$, and in $\S 9$ we will prove Theorem $\ref{t0.4}$.\vspace{5mm}

\noindent \textbf{Remark:} $(\ref{0.1})$ is globally well-posed when $\| u_{0} \|_{L^{2}(\mathbf{R}^{n})}$ is sufficiently small. $\cite{KTV}$ and $\cite{KVZ}$ proved global well-posedness for all $u_{0} \in L^{2}(\mathbf{R}^{n})$, $u_{0}$ radial, based on an induction on mass method. This method will not be used here.

\section{Some Harmonic Analysis}
In this section, the harmonic analysis tools that will be needed later will be given. Let $\mathcal F$ be the Fourier transform,

\begin{equation}\label{1.0}
\mathcal F(f)(\xi) = \int e^{-ix \cdot \xi} f(x) dx.
\end{equation}

\begin{definition}\label{d1.1}
Suppose $\phi(\xi)$ is a $C_{0}^{\infty}$, decreasing, radial function. Also suppose,

\begin{equation}\label{1.1}
\phi(\xi) = \left\{
              \begin{array}{ll}
                1, & \hbox{$|\xi| \leq 1/2$;} \\
                0, & \hbox{$|\xi| > 1$.}
              \end{array}
            \right.
\end{equation}

\noindent Then define the frequency cutoff

\begin{equation}\label{1.2}
\mathcal F(P_{\leq M} f) = \phi(\frac{\xi}{M}) \hat{f}(\xi).
\end{equation}

\begin{equation}\label{1.2.1}
P_{> M} f = f - P_{\leq M} f.
\end{equation}

\begin{equation}\label{1.2.2}
P_{M} f = P_{\leq M} f - P_{\leq \frac{M}{2}} f.
\end{equation}
\end{definition}

\begin{lemma}\label{l1.2}
When $M > N$,

\begin{equation}\label{1.5}
\| P_{> M} u \|_{\dot{H}^{s}(\mathbf{R}^{n})} \lesssim \frac{1}{N^{1 - s}} \| \nabla Iu \|_{L^{2}(\mathbf{R}^{n})}.
\end{equation}
\end{lemma}

\noindent \emph{Proof:} By definition of the I-operator, $$\| \nabla I P_{> 4N} P_{> M} u \|_{L^{2}(\mathbf{R}^{n})} = N^{1 - s} \| |\nabla|^{s} P_{> 4N} P_{> M} u \|_{L^{2}(\mathbf{R}^{n})},$$ so $$\| |\nabla|^{s} P_{> 4N} P_{> M} u \|_{L^{2}(\mathbf{R}^{n})} \lesssim \frac{1}{N^{1 - s}} \| \nabla Iu \|_{L^{2}(\mathbf{R}^{n})}.$$ Meanwhile, $$\| |\nabla|^{s} P_{\leq 4N} P_{> M} u \|_{L^{2}(\mathbf{R}^{n})} \lesssim \sum_{-5 \leq k \leq 5} (N 2^{k})^{s} \| P_{2^{k} N} u \|_{L^{2}(\mathbf{R}^{n})}$$ $$\lesssim \frac{1}{N^{1 - s}} \sum_{-5 \leq k \leq 5} \| \nabla P_{2^{k} N} u \|_{L^{2}(\mathbf{R}^{n})} \lesssim \frac{1}{N^{1 - s}} \| \nabla Iu \|_{L^{2}(\mathbf{R}^{n})}.$$ This proves the lemma. $\Box$

\begin{lemma}\label{l1.3}
Suppose $n \geq 3$. A pair $(p,q)$ is called admissible if

\begin{equation}\label{1.5.1}
\frac{2}{p} = n(\frac{1}{2} - \frac{1}{q}), p \geq 2,
\end{equation}

\begin{equation}\label{1.6}
\| e^{it \Delta} u_{0} \|_{L_{t}^{p} L_{x}^{q}(J \times \mathbf{R}^{n})} \lesssim \| u_{0} \|_{L^{2}(\mathbf{R}^{n})},
\end{equation}

\noindent for pairs $(p,q)$ that satisfy $(\ref{1.5.1})$.
\end{lemma}

\noindent \emph{Proof:} See \cite{Tao} for the case $p > 2$ and \cite{KT} for $p = 2$. $\Box$\vspace{5mm}

\noindent There is a very useful space of functions, called the Strichartz space.

\begin{equation}\label{1.7}
\| u \|_{S^{0}(J \times \mathbf{R}^{n})} = \sup_{(p,q) \text{ admissible}} \| u \|_{L_{t}^{p} L_{x}^{q}(J \times \mathbf{R}^{n})}.
\end{equation}

\noindent Let $p'$ denote the dual exponent to $p$, $\frac{1}{p'} = 1 - \frac{1}{p}$. The dual space to $(\ref{1.7})$ is

\begin{equation}\label{1.8}
\| F \|_{N^{0}(J \times \mathbf{R}^{n})} = \inf_{(p,q) \text{ admissible}} \| F \|_{L_{t}^{p'} L_{x}^{q'}(J \times \mathbf{R}^{n})}.
\end{equation}

\noindent If

\begin{equation}\label{1.8.1}
\aligned
iu_{t} + \Delta u &= F, \\
u(a) &= u_{0},
\endaligned
\end{equation}

\noindent then

\begin{equation}\label{1.9}
\| u \|_{S^{0}([a, b] \times \mathbf{R}^{n})} \lesssim \| u_{0} \|_{L^{2}(\mathbf{R}^{n})} + \| F \|_{N^{0}([a, b] \times \mathbf{R}^{n})}.
\end{equation}

\noindent See \cite{Tao} for more information on this space.\vspace{5mm}

\noindent Finally, the bilinear Strichartz estimate when $n = 3$ will be used to resolve one technical issue.

\begin{theorem}\label{t1.4}
For any spacetime slab $I \times \mathbf{R}^{3}$ and any $t_{0} \in I$, and for any $\delta > 0$, suppose $M << N$, u is supported on frequency $|\xi| \geq N$, v on frequency $|\xi| \leq M$,

\begin{equation}\label{1.10}
\aligned
\| uv \|_{L_{t,x}^{2}(I \times \mathbf{R}^{3})} \leq C(\delta) \frac{M^{1 - \delta}}{N^{1/2 - \delta}} (\| u(t_{0}) \|_{L^{2}(\mathbf{R}^{3})} + \| (i \partial_{t} + \Delta) u \|_{L_{t}^{1} L_{x}^{2}(I \times \mathbf{R}^{3})}) \\
\times (\| v(t_{0}) \|_{L^{2}(\mathbf{R}^{3})} + \| (i \partial_{t} + \Delta) v \|_{L_{t}^{1} L^{2}(I \times \mathbf{R}^{3})})
\endaligned
\end{equation}
\end{theorem}

\noindent \emph{Proof:} See \cite{CKSTT4}.

\section{Smoothness Estimates}
The principle difficulty that arises for $n \geq 3$ is that the nonlinearity $|u|^{4/n} u$ is no longer algebraic. To circumvent this problem, it is necessary to understand how the smoothness of $Iu$ affects the smoothness of $I(|u|^{4/n})$.

\begin{theorem}\label{t2.1}
Suppose $n = 3, 4$. If $M \geq N$,

\begin{equation}\label{2.2}
\aligned
\| P_{> M} |u|^{4/n} \|_{L^{n/2}(\mathbf{R}^{n})} \lesssim \frac{1}{M^{s} N^{1 - s}} \| \langle \nabla \rangle Iu \|_{L^{2}(\mathbf{R}^{n})}^{4/n}, \\
\| P_{> M} \frac{u^{2}}{|u|^{2 - 4/n}} \|_{L^{n/2}(\mathbf{R}^{n})} \lesssim \frac{1}{M^{s} N^{1 - s}} \| \langle \nabla \rangle Iu \|_{L^{2}(\mathbf{R}^{n})}^{4/n}.
\endaligned
\end{equation}

\noindent If $M < N$,

\begin{equation}\label{2.1}
\aligned
\| P_{> M} |u|^{4/n} \|_{L^{n/2}(\mathbf{R}^{n})} \lesssim \frac{1}{M} \| \langle \nabla \rangle Iu \|_{L^{2}(\mathbf{R}^{n})}^{4/n}, \\
\| P_{> M} \frac{u^{2}}{|u|^{2 - 4/n}} \|_{L^{n/2}(\mathbf{R}^{n})} \lesssim \frac{1}{M} \| \langle \nabla \rangle Iu \|_{L^{2}(\mathbf{R}^{n})}^{4/n}.
\endaligned
\end{equation}
\end{theorem}

\noindent \emph{Proof:}

\noindent \emph{Case 1, $M \geq N$:} Split the data, $u = u_{l} + u_{h}$ with $u_{l} = P_{\leq M} u$, $u_{h} = P_{> M} u$. $$\| \nabla u_{l} \|_{L^{2}(\mathbf{R}^{n})} \lesssim \frac{M^{1 - s}}{N^{1 - s}} \| \nabla Iu \|_{L^{2}(\mathbf{R}^{n})}.$$ By elementary calculation and the Leibniz rule

\begin{equation}\label{2.4}
\| \nabla |u_{l}|^{4/n} \|_{L^{n/2}(\mathbf{R}^{n})} \lesssim \frac{M^{1 - s}}{N^{1 - s}} \| \langle \nabla \rangle Iu \|_{L^{2}(\mathbf{R}^{n})}^{4/n},
\end{equation}

\begin{equation}\label{2.5}
\| \nabla \frac{u_{l}^{2}}{|u_{l}|^{2 - 4/n}} \|_{L^{n/2}(\mathbf{R}^{n})} \lesssim \frac{M^{1 - s}}{N^{1 - s}} \| \langle \nabla \rangle Iu \|_{L^{2}(\mathbf{R}^{n})}^{4/n}.
\end{equation}

\noindent Therefore,

\begin{equation}\label{2.6}
\| P_{> M} \frac{u_{l}^{2}}{|u_{l}|^{2 - 4/n}} \|_{L^{n/2}(\mathbf{R}^{n})} + \| P_{> M} |u_{l}|^{4/n} \|_{L^{n/2}(\mathbf{R}^{n})} \lesssim \frac{1}{M} \frac{M^{1 - s}}{N^{1 - s}} \| \langle \nabla \rangle Iu \|_{L^{2}(\mathbf{R}^{n})}^{4/n}.
\end{equation}

\noindent On the other hand, by lemma $\ref{l1.2}$, $M > N$, $$\| u_{h} \|_{L^{2}(\mathbf{R}^{n})} \lesssim \frac{1}{M^{s}} \| |\nabla|^{s} u_{h} \|_{L^{2}(\mathbf{R}^{n})} \lesssim \frac{1}{M^{s}} \frac{1}{N^{1 - s}} \| \langle \nabla \rangle Iu \|_{L^{2}(\mathbf{R}^{n})}.$$ The functions $F(x) = |x|$ and $G(x) = \frac{x^{2}}{|x|}$ are Lipschitz functions, so

\begin{equation}\label{2.7}
\aligned
\| F(u_{l} + u_{h}) - F(u_{l}) \|_{L^{2}(\mathbf{R}^{n})} + \| G(u_{l} + u_{h}) - G(u_{l}) \|_{L^{2}(\mathbf{R}^{n})} \\ \lesssim \| u_{h} \|_{L^{2}(\mathbf{R}^{n})} \lesssim \frac{1}{M^{s} N^{1 - s}} \| \langle \nabla \rangle Iu \|_{L^{2}(\mathbf{R}^{4})}.
\endaligned
\end{equation}

\noindent This takes care of $n = 4$. For $n = 3$ let $F(x) = |x|^{4/3}$ and $G(x) = \frac{x^{2}}{|x|^{2/3}}$.

\begin{equation}\label{2.7.1}
|F(x + y) - F(x)| + |G(x + y) - G(x)| \lesssim |y| (|x|^{1/3} + |y|^{1/3}).
\end{equation}

\noindent Therefore,

\begin{equation}\label{2.8}
\aligned
\| F(u_{l} + u_{h}) - F(u_{l}) \|_{L^{3/2}(\mathbf{R}^{3})} + \| G(u_{l} + u_{h}) - G(u_{l}) \|_{L^{3/2}(\mathbf{R}^{3/2})} \\
\lesssim \| u_{h} \|_{L^{2}(\mathbf{R}^{3})}(\| u_{l} \|_{L^{2}(\mathbf{R}^{3})}^{1/3} + \| u_{h} \|_{L^{2}(\mathbf{R}^{3})}^{1/3}) \lesssim \frac{1}{N^{1 - s} M^{s}} \| \langle \nabla \rangle Iu \|_{L^{2}(\mathbf{R}^{n})}^{4/3}.
\endaligned
\end{equation}

\noindent \emph{Case 2, $M \leq N$:} In this case let $u_{l} = P_{\leq N} u$ and $u_{h} = P_{> N}$. By $(\ref{2.4})$ and $(\ref{2.5})$,

\begin{equation}\label{2.9}
\| P_{> M} \frac{u_{l}^{2}}{|u_{l}|^{2 - 4/n}} \|_{L^{n/2}(\mathbf{R}^{n})} + \| P_{> M} |u_{l}|^{4/n} \|_{L^{n/2}(\mathbf{R}^{n})} \lesssim \frac{1}{M} \| \langle \nabla \rangle Iu \|_{L^{2}(\mathbf{R}^{n})}^{4/n}.
\end{equation}

\noindent Using $(\ref{2.7.1})$ when $n = 3$, $F, G$ Lipschitz when $n = 4$,

\begin{equation}\label{2.10}
\| |u|^{4/n} - |u_{l}|^{4/n} \|_{L^{n/2}(\mathbf{R}^{n})} + \| \frac{u^{2}}{|u|^{2 - 4/n}} - \frac{u_{l}^{2}}{|u_{l}|^{2 - 4/n}} \|_{L^{n/2}(\mathbf{R}^{n})} \lesssim \frac{1}{N} \| \langle \nabla \rangle Iu \|_{L^{2}(\mathbf{R}^{n})}^{4/n}.
\end{equation}

\noindent This completes the proof of the theorem. $\Box$\vspace{5mm}

\noindent When $n > 4$, $u \in H^{1}(\mathbf{R}^{n})$ no longer implies $|u|^{4/n} \in H^{1,p}(\mathbf{R}^{n})$ for any $p$. Instead, it is necessary to rely on a proposition from \cite{V}.

\begin{proposition}\label{p2.2}
Let F be a H{\"o}lder continuous function of order $0 < \alpha < 1$. For every $0 < \sigma < \alpha$, $1 < p < \infty$, $\frac{\sigma}{\alpha} < \rho < 1$,

\begin{equation}\label{2.11}
\| |\nabla|^{\sigma} F(u) \|_{L^{p}(\mathbf{R}^{n})} \lesssim \| |u|^{\alpha - \sigma/\rho} \|_{L^{p_{1}}(\mathbf{R}^{n})} \| |\nabla|^{\rho} u \|_{L^{\frac{\sigma p_{2}}{\rho}}(\mathbf{R}^{n})},
\end{equation}

\noindent provided $\frac{1}{p} = \frac{1}{p_{1}} + \frac{1}{p_{2}}$, $(1 - \frac{\sigma}{\alpha \rho}) p_{1} > 1$.
\end{proposition}

\noindent \emph{Proof:} See \cite{V}.

\begin{theorem}\label{t2.3}
Suppose $n > 4$. If $M > N$,

\begin{equation}\label{2.12}
\| P_{> M} \frac{u^{2}}{|u|^{2 - 4/n}} \|_{L^{n/2}(\mathbf{R}^{n})} + \| P_{> M} |u|^{4/n} \|_{L^{n/2}(\mathbf{R}^{n})} \lesssim \frac{1}{M^{\frac{4s}{n}-}} \frac{1}{N^{\frac{4}{n}(1 - s)-}} \| \langle \nabla \rangle Iu \|_{L^{2}(\mathbf{R}^{n})}^{4/n}.
\end{equation}

\noindent If $M \leq N$,

\begin{equation}\label{2.13}
\| P_{> M} \frac{u^{2}}{|u|^{2 - 4/n}} \|_{L^{2}(\mathbf{R}^{n})} + \| P_{> M} |u|^{4/n} \|_{L^{n/2}(\mathbf{R}^{n})} \lesssim \frac{1}{M^{\frac{4}{n}-}} \| \langle \nabla \rangle Iu \|_{L^{2}(\mathbf{R}^{n})}^{4/n}.
\end{equation}
\end{theorem}

\noindent \emph{Proof:} Let $F(x) = |x|^{4/n}$ and $G(x) = \frac{x^{2}}{|x|^{2 - 4/n}}$, then $F, G \in C^{0, 4/n}(\mathbf{C})$. Choose $\frac{\sigma}{\rho} = \frac{4}{n} - \epsilon$, $\rho = 1 - \delta$.\vspace{5mm}

\noindent \emph{Case 1, $M \geq N$:} Let $u_{l} = P_{\leq M} u$.

\begin{equation}\label{2.14}
\| |\nabla|^{\rho} u_{l} \|_{L^{2}(\mathbf{R}^{n})} \lesssim \frac{M^{1 - s - \delta}}{N^{1 - s - \delta}} \| \langle \nabla \rangle Iu \|_{L^{2}(\mathbf{R}^{n})}.
\end{equation}

\noindent Let $p = \frac{n}{2}$, $p_{1} = \frac{2}{\epsilon}$, $p_{2} = \frac{2n}{4 - n \epsilon}$, $p_{2} \cdot \frac{\sigma}{\rho} = 2$.

\begin{equation}\label{2.15}
\aligned
\| |\nabla|^{\sigma} F(u_{l}) \|_{L^{n/2}(\mathbf{R}^{n})} \lesssim \| |u_{l}|^{\epsilon} \|_{L^{2/\epsilon}(\mathbf{R}^{n})} \| |\nabla|^{\rho} u_{l} \|_{L^{2}(\mathbf{R}^{n})}^{\sigma/\rho} \\
\lesssim \| \langle \nabla \rangle Iu \|_{L^{2}(\mathbf{R}^{n})}^{\epsilon} (\frac{M^{1 - s - \delta}}{N^{1 - s - \delta}})^{\sigma/\rho} \| \langle \nabla \rangle Iu \|_{L^{2}(\mathbf{R}^{n})}^{4/n - \epsilon}.
\endaligned
\end{equation}

\begin{equation}\label{2.16}
\| P_{> M} F(u_{l}) \|_{L^{n/2}(\mathbf{R}^{n})} \lesssim M^{-s(4/n - \epsilon)} N^{(4/n - \epsilon)(1 - s - \delta)} \| \langle \nabla \rangle Iu \|_{L^{2}(\mathbf{R}^{n})}^{4/n}.
\end{equation}

\noindent Make a similar calculation for G. Since F and G are both H{\"o}lder continuous, $$|F(u_{l} + u_{h}) - F(u_{l})| \lesssim |u_{h}|^{4/n},$$ $$|G(u_{l} + u_{h}) - G(u_{l})| \lesssim |u_{h}|^{4/n}.$$ This implies,

\begin{equation}\label{2.17}
\aligned
\| F(u_{l} + u_{h}) - F(u_{l}) \|_{L^{n/2}(\mathbf{R}^{n})} + \| G(u_{l} + u_{h}) - G(u_{l}) \|_{L^{n/2}(\mathbf{R}^{n})} \\ \lesssim (\frac{1}{M^{s} N^{1 - s}})^{4/n} \| \langle \nabla \rangle Iu \|_{L^{2}(\mathbf{R}^{n})}^{4/n}.
\endaligned
\end{equation}

\noindent \emph{Case 2, $M \leq N$:} In this case let $u_{l} = P_{\leq N} u$. In this case

\begin{equation}\label{2.17}
\aligned
\| F(u_{l} + u_{h}) - F(u_{l}) \|_{L^{n/2}(\mathbf{R}^{n})} + \| G(u_{l} + u_{h}) - G(u_{l}) \|_{L^{n/2}(\mathbf{R}^{n})} \\ \lesssim \frac{1}{N^{4/n}} \| \langle \nabla \rangle Iu \|_{L^{2}(\mathbf{R}^{n})}^{4/n}.
\endaligned
\end{equation}

\noindent Also, for $\sigma = (\frac{4}{n} - \epsilon)(1 - \delta)$,

\begin{equation}\label{2.18}
\| |\nabla|^{\sigma} F(u_{l}) \|_{L^{n/2}(\mathbf{R}^{n})} + \| |\nabla|^{\sigma} G(u_{l}) \|_{L^{n/2}(\mathbf{R}^{n})} \lesssim \| \langle \nabla \rangle Iu \|_{L^{2}(\mathbf{R}^{n})}^{4/n}.
\end{equation}

\noindent So in this case,

\begin{equation}\label{2.19}
\| P_{> M} F(u) \|_{L^{n/2}(\mathbf{R}^{n})} + \| P_{> M} G(u) \|_{L^{n/2}(\mathbf{R}^{n})} \lesssim \frac{1}{M^{(4/n - \epsilon)(1 - \delta)}} \| \langle \nabla \rangle Iu \|_{L^{2}(\mathbf{R}^{n})}^{4/n}.
\end{equation}

\noindent Taking $\epsilon, \delta > 0$ arbitrarily small proves the theorem. $\Box$

\section{Linear Term for $n \geq 4$}

\noindent In this section the linear term $(\ref{0.8})$ for $n \geq 4$. The $n = 3$ case is put off until the next section, due to a technical complication.

\begin{theorem}\label{t3.0}
\begin{equation}\label{3.1}
|Im \int_{t_{1}}^{t_{2}} \int (\overline{I \Delta u}) [|Iu|^{4/n} (Iu) - I(|u|^{4/n} u)] dx dt| \lesssim \frac{1}{N^{4/n-}} \| \langle \nabla \rangle Iu \|_{S^{0}(J \times \mathbf{R}^{n})}^{2 + 4/n}.
\end{equation}
\end{theorem}

\noindent To prove this integrate by parts,

\begin{equation}\label{3.1.0}
Im \int_{t_{1}}^{t_{2}} \int (\overline{I \Delta u}) [|Iu|^{4/n} (Iu) - I(|u|^{4/n} u)] dx dt
\end{equation}

\begin{equation}\label{3.2}
= -(\frac{n + 2}{n}) Im \int_{t_{1}}^{t_{2}} \int (\overline{\nabla Iu}) (I(|u|^{4/n} \nabla u) - |Iu|^{4/n} (\nabla Iu)) dx dt
\end{equation}

\begin{equation}\label{3.3}
-\frac{2}{n} Im \int_{t_{1}}^{t_{2}} \int (\overline{\nabla Iu}) (I(\frac{u^{2}}{|u|^{2 - 4/n}} \overline{\nabla u}) - \frac{(Iu)^{2}}{|Iu|^{2 - 4/n}} (\overline{\nabla Iu})) dx dt.
\end{equation}

\noindent First estimate a slightly modified version of $(\ref{3.2})$ and $(\ref{3.3})$, approximating $|Iu|^{4/n}$ by $I(|u|^{4/n})$.

\begin{lemma}\label{l3.1}
Suppose J is an interval.

\begin{equation}\label{3.4}
|\int_{t_{1}}^{t_{2}} \int (\overline{\nabla Iu})(I(|u|^{4/n} \nabla u) - I(|u|^{4/n}) (\nabla Iu)) dx dt| \lesssim \frac{\| \langle \nabla \rangle Iu \|_{S^{0}(J \times \mathbf{R}^{n})}^{2 + 4/n}}{N^{4/n-}}.
\end{equation}

\begin{equation}\label{3.5}
|\int_{t_{1}}^{t_{2}} \int (\overline{\nabla Iu}) (I(\frac{u^{2}}{|u|^{2 - 4/n}} \overline{\nabla u}) - I(\frac{u^{2}}{|u|^{2 - 4/n}}) (\nabla Iu)) dx dt| \lesssim \frac{\| \langle \nabla \rangle Iu \|_{S^{0}(J \times \mathbf{R}^{n})}^{2 + 4/n}}{N^{4/n-}}.
\end{equation}
\end{lemma}

\noindent \emph{Proof:} The proof of $(\ref{3.5})$ is virtually identical to the proof of $(\ref{3.4})$, so only $(\ref{3.4})$ will be proved. Recall that $F(u) = |u|^{4/n}$ and $G(u) = \frac{u^{2}}{|u|^{2 - 4/n}}$. All that will be used to prove $(\ref{3.4})$ is $F \in C^{0, 4/n}(\mathbf{C})$, which is also true for $G$.

$$\int_{t_{1}}^{t_{2}} \int (\overline{\nabla Iu})(I(|u|^{4/n} u) - I(|u|^{4/n}) (\nabla Iu)) dx dt$$

$$ = \int_{t_{1}}^{t_{2}} \int_{\Sigma} (\xi_{1} \widehat{\overline{Iu}}(t, \xi_{1})) (1 - \frac{m(\xi_{2} + \xi_{3})}{m(\xi_{2}) m(\xi_{3})} ) \cdot (\xi_{2} \widehat{Iu}(t,\xi_{2})) (\widehat{IF(u)}(t,\xi_{3})) d\xi dt,$$

\noindent where $\Sigma$ is the hyperplane $\{ \xi_{1} + \xi_{2} + \xi_{3} = 0 \}$ and $d\xi$ is the Lebesgue measure on $\Sigma$. Make a Littlewood - Paley decomposition and consider several cases separately.\vspace{5mm}

\noindent \textbf{Case 1, $N_{2}, N_{3} << N$:} In this case the multiplier is $\equiv 0$.\vspace{5mm}

\noindent \textbf{Case 2(a), $N_{2} \gtrsim N >> N_{3}:$} Here $N_{1} \sim N_{2}$. Using the fundamental theorem of calculus,

$$|1 - \frac{m(\xi_{2} + \xi_{3})}{m(\xi_{2}) m(\xi_{3})}| = |1 - \frac{m(\xi_{2} + \xi_{3})}{m(\xi_{2})}| \lesssim \frac{\nabla m(\xi_{2})}{m(\xi_{2})} |\xi_{3}| \lesssim \frac{N_{3}}{N_{2}}.$$

$$\sum_{N \lesssim N_{1} \sim N_{2}}  \| P_{N_{1}} \nabla Iu \|_{L_{t}^{2} L_{x}^{\frac{2n}{n - 2}}(J \times \mathbf{R}^{n})} \| P_{N_{2}} \nabla Iu \|_{L_{t}^{2} L_{x}^{\frac{2n}{n - 2}}(J \times \mathbf{R}^{n})}$$

$$\times \sum_{N_{3} << N} \frac{N_{3}}{N_{2}} \| P_{N_{3}} IF(u) \|_{L_{t}^{\infty} L_{x}^{n/2}(J \times \mathbf{R}^{n})},$$

\noindent Applying Theorems $\ref{t2.1}$ and $\ref{t2.3}$,

$$\lesssim \sum_{N \lesssim N_{1} \sim N_{2}} \frac{1}{N_{2}} \| P_{N_{1}} \nabla Iu \|_{S^{0}(J \times \mathbf{R}^{n})} \| P_{N_{2}} \nabla Iu \|_{S^{0}(J \times \mathbf{R}^{n})}$$ $$\times \sum_{N_{3} << N} N_{3}^{(1 - 4/n)+} \| \langle \nabla \rangle Iu \|_{S^{0}(J \times \mathbf{R}^{n})}^{4/n}$$

$$\lesssim \frac{1}{N^{4/n-}} \| \langle \nabla \rangle Iu \|_{S^{0}(J \times \mathbf{R}^{n})}^{2 + 4/n}.$$\vspace{5mm}

\noindent \textbf{Case 2(b), $N_{3} \gtrsim N >> N_{2}$:} Making a similar calculation, it suffices to estimate

$$\sum_{N \lesssim N_{1} \sim N_{3}}  \| P_{N_{1}} \nabla Iu \|_{L_{t}^{2} L_{x}^{\frac{2n}{n - 2}}(J \times \mathbf{R}^{n})} \| P_{N_{3}} IF(u) \|_{L_{t}^{\infty} L_{x}^{n/2}(J \times \mathbf{R}^{n})}$$

$$\times \sum_{N_{2} << N} \frac{N_{2}}{N_{3}} \| P_{N_{2}} \nabla Iu \|_{L_{2}^{\infty} L_{x}^{2n/(n - 2)}(J \times \mathbf{R}^{n})},$$

\noindent Again applying Theorems $\ref{t2.1}$ and $\ref{t2.3}$,

$$\lesssim \| \langle \nabla \rangle Iu \|_{S^{0}(J \times \mathbf{R}^{n})}^{2 + 4/n} \sum_{N \lesssim N_{1} \sim N_{3}}  \frac{N}{N_{3}^{4/n-} N_{3}} \lesssim \frac{1}{N^{4/n-}} \| \langle \nabla \rangle Iu \|_{S^{0}(J \times \mathbf{R}^{n})}^{2 + 4/n}. $$\vspace{5mm}

\noindent \textbf{Case 3, $N_{2} \gtrsim N$ and $N_{3} \gtrsim N$}
There are three subcases to consider.\vspace{5mm}

\noindent \textbf{Case 3(a), $N_{1} \sim N_{2}$, $N_{2} >> N_{3}$:} In this case $|1 - \frac{m(\xi_{2} + \xi_{3})}{m(\xi_{2}) m(\xi_{3})}| \lesssim \frac{1}{m(\xi_{3})}$.

$$\sum_{N \lesssim N_{1} \sim N_{2}} \| P_{N_{1}} \nabla Iu \|_{L_{t}^{2} L_{x}^{2n/(n - 2)}(J \times \mathbf{R}^{n})}  \| P_{N_{2}} \nabla Iu \|_{L_{t}^{2} L_{x}^{\frac{2n}{n - 2}}(J \times \mathbf{R}^{n})}$$
$$\times \sum_{N \lesssim N_{3} << N_{2}} \frac{1}{m(N_{3})} \| P_{N_{3}} IF(u) \|_{L_{t}^{\infty} L_{x}^{n/2}(J \times \mathbf{R}^{n})}$$

$$\lesssim \sum_{N \lesssim N_{1} \sim N_{2}}  \| P_{N_{1}} \nabla Iu \|_{S^{0}(J \times \mathbf{R}^{n})}  \| P_{N_{2}} \nabla Iu \|_{S^{0}(J \times \mathbf{R}^{n})}$$ $$\times \sum_{N \lesssim N_{3} << N_{2}} \frac{1}{N^{\frac{4(1 - s)-}{n}} N_{1}^{\frac{4s-}{n}}} \| \langle \nabla \rangle Iu \|_{S^{0}(J \times \mathbf{R}^{n})}^{4/n}$$

$$\lesssim \frac{1}{N^{4/n-}} \| \langle \nabla \rangle Iu \|_{S^{0}(J \times \mathbf{R}^{n})}^{2 + 4/n}.$$

\noindent \textbf{Remark:} The last sum follows by Cauchy Schwartz.\vspace{5mm}

\noindent \textbf{Case 3(b), $N_{1} \sim N_{3}$, $N_{3} >> N_{2}$:} In a similar manner,

$$\sum_{N_{1} \sim N_{3}} \| P_{N_{1}} \nabla Iu \|_{L_{t}^{2} L_{x}^{2n/(n - 2)}(J \times \mathbf{R}^{n})} \|  P_{N_{3}} IF(u) \|_{L_{t}^{\infty} L_{x}^{n/2}(J \times \mathbf{R}^{n})}$$

$$\times \sum_{N \lesssim N_{2} << N_{3}} \frac{1}{m(N_{2})} \| P_{N_{2}} \nabla Iu \|_{L_{t}^{2} L_{x}^{2n/(n - 2)}(J \times \mathbf{R}^{n})}$$

$$\lesssim \| \langle \nabla \rangle Iu \|_{S^{0}(J \times \mathbf{R}^{n})}^{2 + 4/n} \sum_{N \lesssim N_{1} \sim N_{3}} \frac{1}{N_{3}^{4/n-}} \frac{N_{3}^{1 - s}}{N^{1 - s}}$$

$$\lesssim \frac{1}{N^{4/n-}} \| \langle \nabla \rangle Iu \|_{S^{0}(J \times \mathbf{R}^{n})}^{2 + 4/n}.$$

\noindent \textbf{Remark:} This proof utilizes the fact that $1 - s < \frac{4}{n}$ in Theorem $\ref{t2.3}$ .\vspace{5mm}

\noindent \textbf{Case 3(c), $N_{2} \sim N_{3}$, $N_{2} \gtrsim N_{1}$:} In this case $|1 - \frac{m(\xi_{2} + \xi_{3})}{m(\xi_{2}) m(\xi_{3})}| \lesssim \frac{1}{m(\xi_{2}) m(\xi_{3})}$.

$$\sum_{N \lesssim N_{2} \sim N_{3}} \| P_{N_{2}} \nabla Iu \|_{L_{t}^{2} L_{x}^{2n/(n - 2)}(J \times \mathbf{R}^{n})} \|  P_{N_{3}} IF(u) \|_{L_{t}^{\infty} L_{x}^{n/2}(J \times \mathbf{R}^{n})}$$

$$\times \sum_{N \lesssim N_{1} \lesssim N_{2}} \frac{m(N_{1})}{m(N_{2})^{2}} \| P_{N_{1}} \nabla Iu \|_{L_{t}^{2} L_{x}^{2n/(n - 2)}(J \times \mathbf{R}^{n})}$$

$$\lesssim \|  \langle \nabla \rangle Iu \|_{S^{0}(J \times \mathbf{R}^{n})}^{2 + 4/n} \sum_{N \lesssim N_{2}} \frac{\ln(N)}{m(N_{2})^{2} N_{3}^{4/n-}}$$

$$\lesssim \frac{1}{N^{4/n-}} \| \langle \nabla \rangle Iu \|_{S^{0}(J \times \mathbf{R}^{n})}^{2 + 4/n}.$$

\noindent When $n \geq 4$, $2(1 - s) < \frac{4}{n}$ since $s > \frac{n - 2}{n}$. This takes care of the lemma. $\Box$\vspace{5mm}

\noindent \textbf{Remark:} It is in this particular case where the above argument would break down when $n = 3$. Therefore, $n = 3$ requires a different method.\vspace{5mm}

\noindent To finish the proof of theorem $\ref{t3.0}$, it remains to prove $$IF(u) + IG(u)$$ is a good approximation of $$F(Iu) + G(Iu).$$

$$(\overline{\nabla Iu}) I(|u|^{4/n} \nabla u) - |Iu|^{4/n} (\nabla Iu) (\overline{\nabla Iu})$$

$$ = (\overline{\nabla Iu}) I(|u|^{4/n} \nabla u) - (\overline{\nabla Iu})(\nabla Iu) I(|u|^{4/n})$$ $$+ (\overline{\nabla Iu})(\nabla Iu) I(|u|^{4/n}) - |Iu|^{4/n} (\nabla Iu) (\overline{\nabla Iu}),$$

\noindent and similarly for $G(u)$.

\begin{lemma}\label{l3.2}
\begin{equation}\label{3.6}
|\int_{t_{1}}^{t_{2}} \int |\nabla Iu|^{2} [I(|u|^{4/n}) - |Iu|^{4/n}] dx dt| \lesssim \frac{1}{N^{4/n-}} \| \langle \nabla \rangle Iu \|_{S^{0}(J \times \mathbf{R}^{n})}^{2 + 4/n},
\end{equation}

\begin{equation}\label{3.8}
|\int_{t_{1}}^{t_{2}} \int (\overline{\nabla Iu})^{2} [I(\frac{u^{2}}{|u|^{2 - 4/n}}) - \frac{(Iu)^{2}}{|Iu|^{2 - 4/n}} dx dt| \lesssim \frac{1}{N^{4/n-}} \| \langle \nabla \rangle Iu \|_{S^{0}(J \times \mathbf{R}^{n})}^{2 + 4/n}.
\end{equation}
\end{lemma}

\noindent \emph{Proof:} Split the data $u = u_{l} + u_{h}$, with $u_{l} = P_{\leq \frac{N}{4}} u$, in particular $Iu_{l} = u_{l}$, and

\begin{equation}\label{3.1}
\| |u_{h}|^{4/n} \|_{L_{t}^{\infty} L_{x}^{n/2}(J \times \mathbf{R}^{n})} \lesssim \frac{1}{N^{4/n}} \| \langle \nabla \rangle Iu \|_{S^{0}(J \times \mathbf{R}^{n})}.
\end{equation}

$$IF(u) - F(Iu) = [IF(u) - IF(u_{l})] + [IF(u_{l}) - F(Iu_{l})] + [F(Iu_{l}) - F(Iu)],$$

$$IG(u) - G(Iu) = [IG(u) - IG(u_{l})] + [IG(u_{l}) - G(Iu_{l})] + [G(Iu_{l}) - G(Iu)].$$

\noindent Since $F, G \in C^{0, 4/n}$,

\begin{equation}\label{3.9}
\aligned
\| F(u) - F(u_{l}) \|_{L_{t}^{\infty} L_{x}^{n/2}(J \times \mathbf{R}^{n})} + \| G(u) - G(u_{l}) \|_{L_{t}^{\infty} L_{x}^{n/2}(J \times \mathbf{R}^{n})} \\ \lesssim \| u_{h} \|_{L_{t}^{\infty} L_{x}^{2}(J \times \mathbf{R}^{n})}^{4/n} \lesssim \frac{1}{N^{4/n}} \| \langle \nabla \rangle Iu \|_{S^{0}(J \times \mathbf{R}^{n})}.
\endaligned
\end{equation}

\noindent Similarly, since $Iu_{l} = u_{l}$;

\begin{equation}\label{3.9.1}
\aligned
\| F(Iu) - F(u_{l}) \|_{L_{t}^{\infty} L_{x}^{n/2}(J \times \mathbf{R}^{n})} + \| G(Iu) - G(u_{l}) \|_{L_{t}^{\infty} L_{x}^{n/2}(J \times \mathbf{R}^{n})} \\
\lesssim \frac{1}{N^{4/n}} \| \langle \nabla \rangle Iu \|_{S^{0}(J \times \mathbf{R}^{n})},
\endaligned
\end{equation}

\noindent Finally,

\begin{equation}\label{3.10}
\aligned
&\| I(|u_{l}|^{4/n}) - |u_{l}|^{4/n} \|_{L_{t}^{\infty} L_{x}^{n/2}(J \times \mathbf{R}^{n})} + \| I(\frac{u_{l}^{2}}{|u_{l}|^{2 - 4/n}}) - \frac{u_{l}^{2}}{|u_{l}|^{2 - 4/n}} \|_{L_{t}^{\infty} L_{x}^{n/2}(J \times \mathbf{R}^{n})} \\ &\lesssim \frac{1}{N^{4/n-}} \| \langle \nabla \rangle Iu \|_{S^{0}(J \times \mathbf{R}^{n})},
\endaligned
\end{equation}

\noindent by $m(\xi) \equiv 1$ on $|\xi| \leq N$, theorems $\ref{t2.1}$ and $\ref{t2.3}$. This proves the lemma. $\Box$\vspace{5mm}

\noindent Combining Lemma $\ref{l3.1}$ and Lemma $\ref{l3.2}$ proves Theorem $\ref{t3.0}$ $\Box$.

\section{Linear Term for $n = 3$}
When $n = 3$, it is necessary to use a different method than was used for $n \geq 4$.

\begin{theorem}\label{t4.1}
\begin{equation}\label{4.1}
|Im \int_{t_{1}}^{t_{2}} \int (\overline{I \Delta u}) [|Iu|^{4/3} (Iu) - I(|u|^{4/3} u)] dx dt| \lesssim \frac{1}{N^{1-}} \| \langle \nabla \rangle Iu \|_{S^{0}(J \times \mathbf{R}^{3})}^{7/3}.
\end{equation}
\end{theorem}

\noindent \emph{Proof:} Let $u = u_{l} + u_{h}$, $u_{l} = P_{\leq N/4} u$. Then $Iu_{l} = u_{l}$. Integrating by parts,

\begin{equation}\label{4.2}
\aligned
\int_{t_{1}}^{t_{2}} \int (\overline{I \Delta u}) &[|u_{l}|^{4/n} (u_{l}) - I(|u_{l}|^{4/n} u_{l})] dx dt, \\
 = \int_{t_{1}}^{t_{2}} \int (\overline{\nabla Iu}) [&I(\frac{5}{3} |u_{l}|^{4/3} \nabla u_{l} + \frac{2}{3} \frac{u_{l}^{2}}{|u_{l}|^{2/3}} \nabla \bar{u_{l}}) \\ &- (\frac{5}{3} |u_{l}|^{4/3} \nabla u_{l} + \frac{2}{3} \frac{u_{l}^{2}}{|u_{l}|^{2/3}} \nabla \bar{u_{l}})] dx dt.
 \endaligned
\end{equation}

\noindent When $N_{1} \leq N$, $$P_{N_{1}}[I(\frac{5}{3} |u_{l}|^{4/3} \nabla u_{l} + \frac{2}{3} \frac{u_{l}^{2}}{|u_{l}|^{2/3}} \nabla \bar{u_{l}}) - (\frac{5}{3} |u_{l}|^{4/3} \nabla u_{l} + \frac{2}{3} \frac{u_{l}^{2}}{|u_{l}|^{2/3}} \nabla \bar{u_{l}})] \equiv 0.$$

\noindent For $N_{1} \geq N$, since $\nabla u_{l}$ is supported on $|\xi| \leq \frac{N}{4}$ it suffices to estimate $$\frac{5}{3} \| P_{|\xi| \sim N_{1}} |u_{l}|^{4/3} \|_{L_{t}^{\infty} L_{x}^{3/2}(J \times \mathbf{R}^{n})} + \frac{2}{3} \| P_{|\xi| \sim N_{1}} \frac{u_{l}^{2}}{|u_{l}|^{2/3}} \|_{L_{t}^{\infty} L_{x}^{3/2}(J \times \mathbf{R}^{n})}$$

$$\lesssim \frac{1}{N_{1}} \| \nabla P_{|\xi| \sim N_{1}} |u_{l}|^{4/3} \|_{L_{t}^{\infty} L_{x}^{3/2}(J \times \mathbf{R}^{n})} + \frac{1}{N_{1}} \| \nabla P_{|\xi| \sim N_{1}} \frac{u_{l}^{2}}{|u_{l}|^{2/3}} \|_{L_{t}^{\infty} L_{x}^{3/2}(J \times \mathbf{R}^{n})}$$

$$\lesssim \frac{1}{N_{1}} \| \langle \nabla \rangle Iu \|_{S^{0}(J \times \mathbf{R}^{n})}^{4/3},$$ so,

\begin{equation}\label{4.3}
(\ref{4.2}) \lesssim \sum_{N \lesssim N_{1}} \frac{1}{N_{1}} \| \langle \nabla \rangle Iu \|_{S^{0}(J \times \mathbf{R}^{n})}^{10/3} \lesssim \frac{1}{N} \| \langle \nabla \rangle Iu \|_{S^{0}(J \times \mathbf{R}^{n})}^{10/3}.
\end{equation}

\noindent Next, use the Taylor expansion, $$f(x + y) = f(x) + \int_{0}^{1} y f'(x + \tau y) d\tau.$$

\begin{equation}\label{4.4}
|Iu|^{4/3} (Iu) = |u_{l}|^{4/3} u_{l} + \int_{0}^{1} \frac{5}{3} |u_{l} + \tau Iu_{h}|^{4/3} (Iu_{h}) + \frac{2}{3} \frac{(u_{l} + \tau Iu_{h})^{2}}{|u_{l} + \tau Iu_{h}|^{2/3}} (\overline{Iu_{h}}) d\tau.
\end{equation}

\begin{equation}\label{4.5}
I(|u|^{4/3} u) = I(|u_{l}|^{4/3} u_{l}) + \int_{0}^{1} \frac{5}{3} I(|u_{l} + \tau u_{h}|^{4/3} (u_{h})) + \frac{2}{3} I(\frac{(u_{l} + \tau u_{h})^{2}}{|u_{l} + \tau u_{h}|^{2/3}} (\overline{u_{h}})) d\tau.
\end{equation}

$$\int_{t_{1}}^{t_{2}} \int (\Delta Iu) [I(F(u_{l} + \tau u_{h}) u_{h}) - F(u_{l} + \tau Iu_{h}) Iu_{h}] dx dt$$

\begin{equation}\label{4.6}
\aligned
 = \int_{t_{1}}^{t_{2}} \int_{\Sigma} (|\xi_{1}|^{2} \widehat{Iu}(t,\xi_{1}))  [m(\xi_{2} + \xi_{3}) \hat{F}(u_{l} + \tau u_{h})(t, \xi_{2}) \hat{u}_{h}(t,\xi_{3}) \\ - m(\xi_{3}) \hat{F}(u_{l} + \tau Iu_{h})(t, \xi_{2}) (\widehat{Iu_{h}})(t,\xi_{3})] d\xi dt.
\endaligned
\end{equation}

\noindent As usual, make a Littlewood - Paley decomposition and consider several cases separately. It suffices to consider only $N_{3} \gtrsim N$ because of the support of $u_{h}$.\vspace{5mm}

\noindent \textbf{Case 1, $N_{1} \sim N_{3} \gtrsim N$, $N_{2} << N$:}

\begin{equation}\label{4.7}
\aligned
(\ref{4.6}) = \int_{t_{1}}^{t_{2}} \int_{\Sigma} (|\xi_{1}|^{2} \widehat{Iu}(t,\xi_{1}))  [m(\xi_{2} + \xi_{3}) \hat{F}(u_{l} + \tau u_{h})(t, \xi_{2}) \hat{u}_{h}(t,\xi_{3}) \\ - m(\xi_{3}) \hat{F}(u_{l} + \tau u_{h})(t, \xi_{2}) (\widehat{Iu_{h}})(t,\xi_{3})] d\xi dt
\endaligned
\end{equation}

\begin{equation}\label{4.8}
\aligned
+ \int_{t_{1}}^{t_{2}} \int_{\Sigma} (|\xi_{1}|^{2} \widehat{Iu}(t,\xi_{1})) m(\xi_{3})[\hat{F}(u_{l} + \tau u_{h})(t, \xi_{2}) \\- \hat{F}(u_{l} + \tau Iu_{h})(t, \xi_{2})] (\widehat{Iu_{h}})(t,\xi_{3})] d\xi dt.
\endaligned
\end{equation}

\noindent For $(\ref{4.7})$, using the fundamental theorem of calculus, $$|m(N_{2} + N_{3}) - m(N_{3})| \lesssim \frac{N_{2} m(N_{3})}{N_{3}} .$$

$$\sum_{N \lesssim N_{1} \sim N_{3}} \| P_{N_{1}} \Delta Iu \|_{L_{t}^{2} L_{x}^{6}(J \times \mathbf{R}^{3})} \frac{m(N_{3})}{N_{3}} \| P_{N_{3}} u_{h} \|_{L_{t}^{2} L_{x}^{6}(J \times \mathbf{R}^{3})}$$ $$\times \sum_{N_{2} << N} N_{2} \| P_{N_{2}} F(u_{l} + \tau u_{h}) \|_{L_{t}^{\infty} L_{x}^{3/2}(J \times \mathbf{R}^{3})}$$

$$\lesssim \| \langle \nabla \rangle Iu \|_{S^{0}(J \times \mathbf{R}^{3})}^{10/3} \sum_{N \lesssim N_{1} \sim N_{3}} \frac{N_{1} ln(N)}{N_{3}^{2}} \lesssim \frac{1}{N^{1-}} \| \langle \nabla \rangle Iu \|_{S^{0}(J \times \mathbf{R}^{3})}^{10/3},$$

\noindent by theorem $\ref{t2.1}$. To estimate $(\ref{4.8})$, $$|F(u_{l} + \tau u_{h}) - F(u_{l} + \tau Iu_{h})| \lesssim |(1 - I)u_{h}| (|u_{l}|^{1/3} + |u_{h}|^{1/3}),$$

\noindent therefore,

\begin{equation}\label{4.9}
\| P_{N_{2}} [F(u_{l} + \tau u_{h}) - F(u_{l} + \tau Iu_{h})] \|_{L_{t}^{\infty} L_{x}^{3/2}(J \times \mathbf{R}^{3})} \lesssim \frac{1}{N} \| \langle \nabla \rangle Iu \|_{S^{0}(J \times \mathbf{R}^{3})}^{4/3}.
\end{equation}

$$(\ref{4.8}) \lesssim \frac{\ln(N)}{N} \| \langle \nabla \rangle Iu \|_{S^{0}(J \times \mathbf{R}^{3})}^{4/3} \sum_{N \lesssim N_{1} \sim N_{3}} \frac{N_{1}}{N_{3}} \| P_{N_{1}} \nabla Iu \|_{L_{t}^{2} L_{x}^{6}(J \times \mathbf{R}^{3})} \| P_{N_{3}} \nabla Iu \|_{L_{t}^{2} L_{x}^{6}(J \times \mathbf{R}^{3})}$$

$$\lesssim \frac{1}{N^{1-}} \| \langle \nabla \rangle Iu \|_{S^{0}(J \times \mathbf{R}^{3})}^{10/3}.$$

\noindent \textbf{Remark:} Summing in $N_{1} \sim N_{3}$ follows by Cauchy-Schwartz.

\noindent \textbf{Case 2: $N_{1}, N_{2}, N_{3} \gtrsim N$} In this case consider $I(F(u) u_{h})$ and $F(Iu)(Iu_{h})$ separately.\vspace{5mm}

\noindent \textbf{Case 2(a): $N_{1} \sim N_{3} >> N_{2} \gtrsim N$} For the $I(F(u) u_{h})$ term, $m(\xi_{2} + \xi_{3}) \sim m(\xi_{3})$. Using theorem $\ref{t2.1}$ again,

$$\sum_{N \lesssim N_{1} \sim N_{3}} \| P_{N_{1}} \Delta Iu \|_{L_{t}^{2} L_{x}^{6}(J \times \mathbf{R}^{3})} \| P_{N_{3}} Iu \|_{L_{t}^{2} L_{x}^{6}(J \times \mathbf{R}^{3})} \sum_{N \lesssim N_{2} << N_{3}} \| P_{N_{2}} F(u) \|_{L_{t}^{\infty} L_{x}^{3/2}(J \times \mathbf{R}^{3})}$$

$$\lesssim \| \langle \nabla \rangle Iu \|_{S^{0}(J \times \mathbf{R}^{3})}^{4/3} \sum_{N \lesssim N_{1} \sim N_{3}} \frac{N_{1}}{N_{3}} \| P_{N_{1}} \nabla Iu \|_{L_{t}^{2} L_{x}^{6}(J \times \mathbf{R}^{3})}$$ $$\times \| P_{N_{3}} \nabla Iu \|_{L_{t}^{2} L_{x}^{6}(J \times \mathbf{R}^{3})}  \sum_{N \lesssim N_{2} << N_{3}} \frac{1}{N^{1 - s} N_{2}^{s}}$$

$$\lesssim \frac{1}{N^{1-}} \| \langle \nabla \rangle Iu \|_{S^{0}(J \times \mathbf{R}^{3})}.$$

\noindent For the $F(Iu)(Iu)$ term,

$$\sum_{N \lesssim N_{1} \sim N_{3}} \| P_{N_{1}} \Delta Iu \|_{L_{t}^{2} L_{x}^{6}(J \times \mathbf{R}^{3})} \| P_{N_{3}} Iu \|_{L_{t}^{2} L_{x}^{6}(J \times \mathbf{R}^{3})} \sum_{N \lesssim N_{2} << N_{3}} \| P_{N_{3}} F(Iu) \|_{L_{t}^{\infty} L_{x}^{3/2}(J \times \mathbf{R}^{3})}$$

$$\lesssim \| \langle \nabla \rangle Iu \|_{S^{0}(J \times \mathbf{R}^{3})}^{4/3} \sum_{N \lesssim N_{1} \sim N_{3}} \frac{N_{1}}{N_{3}} \| P_{N_{1}} \nabla Iu \|_{L_{t}^{2} L_{x}^{6}(J \times \mathbf{R}^{3})}$$ $$\times \| P_{N_{3}} \nabla Iu \|_{L_{t}^{2} L_{x}^{6}(J \times \mathbf{R}^{3})}  \sum_{N \lesssim N_{2} << N_{3}} \frac{1}{N_{2}}$$

$$\lesssim \frac{1}{N^{1-}} \| \langle \nabla \rangle Iu \|_{S^{0}(J \times \mathbf{R}^{3})}.$$

\noindent \textbf{Case 2(b): $N_{1} \sim N_{2} >> N_{3} \gtrsim N$} In this case $m(\xi_{2} + \xi_{3}) \sim m(\xi_{2})$. To estimate

\begin{equation}\label{4.10}
\aligned
\int_{t_{1}}^{t_{2}} \int \sum_{N \lesssim N_{1} \sim N_{2}}  (P_{N_{1}} \Delta Iu) ( P_{N_{2}} I F(u) ) \sum_{N \lesssim N_{3} << N_{1}} ( P_{N_{3}} u_{h} ),
\endaligned
\end{equation}

\noindent there is a slight technical complication due to the fact that Cauchy - Schwartz is not available for $P_{N_{2}} IF(u)$ ($\S 3$ only proved an estimate on the decay of $P_{N_{2}} IF(u)$, it did not prove $IF(u) \in H^{1, 3/2}(\mathbf{R}^{3})$). Therefore, it is necessary to utilize the bilinear estimates of theorem $\ref{t1.4}$. Interpolating

$$(\ref{4.10}) \lesssim \sum_{N \lesssim N_{1} \sim N_{2}} \sum_{N \lesssim N_{3} << N} \| (P_{N_{1}} \Delta Iu) (P_{N_{3}} u_{h}) \|_{L_{t}^{4/3} L_{x}^{2}(J \times \mathbf{R}^{3})} \| P_{N_{2}} IF(u) \|_{L_{t}^{4} L_{x}^{2}(J \times \mathbf{R}^{3})}.$$

\noindent with the bilinear Strichartz estimate

\begin{equation}\label{4.11}
\aligned
\| (P_{N_{1}} \Delta Iu) (P_{N_{3}} u_{h}) \|_{L_{t,x}^{2}(J \times \mathbf{R}^{3})} \leq C(\delta) \frac{N_{3}^{1 - \delta}}{N_{1}^{1/2 - \delta}} \frac{N_{1}}{N_{3}^{s} N^{1 - s}} \\ \times (\| \langle \nabla \rangle Iu \|_{S^{0}(J \times \mathbf{R}^{3})}^{2} + \| \langle \nabla \rangle Iu \|_{S^{0}(J \times \mathbf{R}^{3})}^{14/3}).
\endaligned
\end{equation}

\noindent Let $\frac{1}{p} = \frac{3 - \epsilon}{4}$, $$\| (P_{N_{1}} \Delta Iu) (P_{N_{3}} u_{h}) \|_{L_{t}^{p} L_{x}^{2}(J \times \mathbf{R}^{3})}$$ $$\lesssim \frac{N_{3}^{\epsilon - \delta \epsilon}}{N_{1}^{\epsilon/2 - \epsilon \delta}} \frac{N_{1}}{N_{3}^{s} N^{1 - s}} (\| \langle \nabla \rangle Iu \|_{S^{0}(J \times \mathbf{R}^{3})}^{2} + \| \langle \nabla \rangle Iu \|_{S^{0}(J \times \mathbf{R}^{3})}^{2 + 8 \epsilon/3}).$$

\noindent Suppose $|J| \lesssim N^{\alpha}$ for some $\alpha$,

$$\| (P_{N_{1}} \Delta Iu)(P_{N_{3}} u_{h}) \|_{L_{t}^{4/3} L_{x}^{2}(J \times \mathbf{R}^{3})}$$ $$\lesssim N^{\epsilon \alpha/4} \frac{N_{3}^{\epsilon - \epsilon \delta}}{N_{1}^{\epsilon/2 - \epsilon \delta}} \frac{N_{1}}{N_{3}^{s} N^{1 - s}} (\| \langle \nabla \rangle Iu \|_{S^{0}(J \times \mathbf{R}^{3})}^{2} + \| \langle \nabla \rangle Iu \|_{S^{0}(J \times \mathbf{R}^{3})}^{2 + 14 \epsilon/3}).$$

$$\| \langle \nabla \rangle Iu \|_{S^{0}(J \times \mathbf{R}^{3})}^{10/3} N^{\epsilon \alpha/4} \sum_{N \lesssim N_{1} \sim N_{2}} \frac{N_{1}}{N_{2} N_{1}^{\epsilon/2 - \epsilon \delta}} \sum_{N \lesssim N_{3} << N_{2}} \frac{N_{3}^{\epsilon - \epsilon \delta}}{N_{3}^{s} N^{1 - s}}$$

$$\lesssim \| \langle \nabla \rangle Iu \|_{S^{0}(J \times \mathbf{R}^{3})}^{10/3} \frac{N^{\epsilon \alpha/4} N^{\epsilon/2}}{N}.$$

\noindent Letting $\epsilon \searrow 0$ proves the claim.\vspace{5mm}

\noindent \textbf{Remark:} In $\S 7$ we will rescale to make $\| \langle \nabla \rangle Iu \|_{S^{0}(J \times \mathbf{R}^{3})} \lesssim 1$, so it will not be necessary to worry about the $ \| \langle \nabla \rangle Iu \|_{S^{0}(J \times \mathbf{R}^{3})}^{2 + 8 \epsilon/3}$ term here, since it will be $\lesssim \| \langle \nabla \rangle Iu \|_{S^{0}(J \times \mathbf{R}^{3})}^{2}$. This rescaling will rescale the interval $[0, T_{0}]$ to $[0, N^{\frac{2(1 - s)}{s}} T_{0}]$, so $|J| \lesssim N^{\alpha}$.\vspace{5mm}

\noindent For the $F(Iu) (Iu_{h})$ term, $$\| \nabla F(Iu) \|_{L_{t}^{\infty} L_{x}^{3/2}(J \times \mathbf{R}^{3})} \lesssim \| \nabla Iu \|_{L_{t}^{\infty} L_{x}^{2}(J \times \mathbf{R}^{3})}^{4/3}.$$

$$\sum_{N \lesssim N_{1} \sim N_{2}} \| P_{N_{1}} \Delta Iu \|_{L_{t}^{2} L_{x}^{6}(J \times \mathbf{R}^{3})} \| P_{N_{2}} F(Iu) \|_{L_{t}^{\infty} L_{x}^{3/2}}$$

$$\times \sum_{N \lesssim N_{3} << N_{1}} \| P_{N_{3}} Iu_{h} \|_{L_{t}^{2} L_{x}^{6}(J \times \mathbf{R}^{3})}$$

$$\lesssim \frac{1}{N} \| \langle \nabla \rangle Iu \|_{S^{0}(J \times \mathbf{R}^{3})} \sum_{N \lesssim N_{1} \sim N_{2}} \frac{N_{1}}{N_{2}} \| P_{N_{1}} \nabla Iu \|_{L_{t}^{2} L_{x}^{6}(J \times \mathbf{R}^{3})} \| P_{N_{2}} \nabla F(Iu) \|_{L_{t}^{\infty} L_{x}^{3/2}}$$

$$\lesssim \frac{1}{N^{1-}} \| \langle \nabla \rangle Iu \|_{S^{0}(J \times \mathbf{R}^{3})}^{10/3}.$$

\noindent In this case $F(Iu) \in H^{1,3/2}(\mathbf{R}^{3})$, so it is possible to use Cauchy - Schwartz.\vspace{5mm}

\noindent \textbf{Case 2(c), $N_{2} \sim N_{3} \gtrsim N_{1} \gtrsim N$:} For the $I(F(u) u_{h})$ term, use the fact that $m(\xi_{2} + \xi_{3}) = m(\xi_{1})$.

$$\sum_{N \lesssim N_{2} \sim N_{3}} \| P_{N_{3}} u_{h} \|_{L_{t}^{2} L_{x}^{6}(J \times \mathbf{R}^{3})} \| P_{N_{2}} F(u) \|_{L_{t}^{\infty} L_{x}^{3/2}(J \times \mathbf{R}^{3})} \sum_{N \lesssim N_{1} \lesssim N_{2}} \| P_{N_{1}} \Delta I^{2} u \|_{L_{t}^{2} L_{x}^{6}(J \times \mathbf{R}^{3})}$$

$$\lesssim \| \langle \nabla \rangle Iu \|_{S^{0}(J \times \mathbf{R}^{3})}^{10/3} \sum_{N \lesssim N_{2} \sim N_{3}} \frac{1}{N_{3}^{s} N_{2}^{s} N^{2(1 - s)}} \sum_{N \lesssim N_{1} \lesssim N_{2}} N_{1}^{s} N^{1 - s}$$

$$\lesssim \| \langle \nabla \rangle Iu \|_{S^{0}(J \times \mathbf{R}^{3})}^{10/3} \sum_{N \lesssim N_{2}} \frac{1}{N_{2}^{s} N^{1 - s}} \lesssim \frac{1}{N^{1-}} \| \langle \nabla \rangle Iu \|_{S^{0}(J \times \mathbf{R}^{3})}^{10/3}.$$\vspace{5mm}

\noindent For the $F(Iu)(Iu)$ term,

$$\sum_{N \lesssim N_{2} \sim N_{3}} \| P_{N_{3}} Iu_{h} \|_{L_{t}^{2} L_{x}^{6}(J \times \mathbf{R}^{3})} \| P_{N_{2}} F(Iu) \|_{L_{t}^{\infty} L_{x}^{3/2}(J \times \mathbf{R}^{3})} \sum_{N \lesssim N_{1} \lesssim N_{2}} \| P_{N_{1}} \Delta Iu \|_{L_{t}^{2} L_{x}^{6}(J \times \mathbf{R}^{3})}$$

$$\lesssim \| \langle \nabla \rangle Iu \|_{S^{0}(J \times \mathbf{R}^{3})}^{7/3} \sum_{N \lesssim N_{2} \sim N_{3}} \frac{1}{N_{2} N_{3}} \sum_{N \lesssim N_{1} \lesssim N_{2}} N_{1}$$

$$\lesssim \frac{1}{N^{1-}} \| \langle \nabla \rangle Iu \|_{S^{0}(J \times \mathbf{R}^{3})}^{7/3}.$$\vspace{5mm}

\noindent A similar calculation can be made for the G term. This proves Theorem $\ref{t4.1}$. $\Box$

\section{Nonlinear Estimate}

\noindent Having dealt with $(\ref{0.7})$, we turn our attention to $(\ref{0.8})$.

\begin{theorem}\label{t5.1}
\begin{equation}\label{5.1}
\aligned
|\int_{t_{1}}^{t_{2}} \int I(\overline{|u|^{4/n} u})(I(|u|^{4/n} u) - |Iu|^{4/n}(Iu)) dx dt| \lesssim \\
\left\{
  \begin{array}{ll}
    \frac{1}{N^{1-}} \| \langle \nabla \rangle Iu \|_{S^{0}(J \times \mathbf{R}^{3})}^{14/3}, & \hbox{when $n = 3$;} \\
    \frac{1}{N^{\frac{4}{n}-}} \| \langle \nabla \rangle Iu \|_{S^{0}(J \times \mathbf{R}^{3})}^{2 + 8/n}, & \hbox{when $n \geq 4$.}
  \end{array}
\right.
\endaligned
\end{equation}
\end{theorem}

\noindent \emph{Proof:} For simplicity of notation let $\sigma(n)$ be the exponent for $N$ in $(\ref{5.1})$. It suffices to prove

\begin{equation}\label{5.3}
\| I(|u|^{4/n} u) - |Iu|^{4/n} (Iu) \|_{L_{t,x}^{2}(J \times \mathbf{R}^{n})} \lesssim \frac{1}{N^{\sigma(n)}} \| \langle \nabla \rangle Iu \|_{S^{0}(J \times \mathbf{R}^{n})}^{1 + 4/n},
\end{equation}

\noindent since $$\| I(|u|^{4/n} u) \|_{L_{t,x}^{2}(J \times \mathbf{R}^{n})} \lesssim \| \nabla I(|u|^{4/n} u) \|_{L_{t}^{2} L_{x}^{\frac{2n}{n + 2}}(J \times \mathbf{R}^{n})}$$

$$\lesssim \| \nabla Iu \|_{L_{t}^{2} L_{x}^{\frac{2n}{n - 2}}(J \times \mathbf{R}^{n})} \| |u|^{4/n} \|_{L_{t}^{\infty} L_{x}^{n/2}(J \times \mathbf{R}^{n})} \lesssim \| \langle \nabla \rangle Iu \|_{S^{0}(J \times \mathbf{R}^{n})}^{1 + 4/n}.$$

\noindent As in the linear case, the quantity in $(\ref{5.3})$ will be split into a main term and a remainder term. This time, we will deal with the remainder term first.

\begin{lemma}\label{l5.2}
\begin{equation}\label{5.4}
\| [|Iu|^{4/n} - |u|^{4/n}](Iu) \|_{L_{t,x}^{2}(J \times \mathbf{R}^{n})} \lesssim \frac{1}{N^{\sigma}} \| \langle \nabla \rangle Iu \|_{S^{0}(J \times \mathbf{R}^{n})}^{1 + 4/n}.
\end{equation}
\end{lemma}

\noindent \emph{Proof:} First consider $n > 4$. $$\| [|Iu|^{4/n} - |u|^{4/n}] (Iu) \|_{L_{t,x}^{2}(J \times \mathbf{R}^{n})} \lesssim \| Iu \|_{L_{t}^{2} L_{x}^{\frac{2n}{n - 4}}(J \times \mathbf{R}^{n})} \| |Iu|^{4/n} - |u|^{4/n} \|_{L_{t}^{\infty} L_{x}^{n/2}(J \times \mathbf{R}^{n})}$$

$$\lesssim \| \langle \nabla \rangle Iu \|_{S^{0}(J \times \mathbf{R}^{n})} \| |Iu|^{4/n} - |u|^{4/n} \|_{L_{t}^{\infty} L_{x}^{n/2}(J \times \mathbf{R}^{n})},$$

\noindent Let $u_{l} = P_{\leq N} u$, since $F(x) \in C^{0, 4/n}$,

$$\lesssim \| \langle \nabla \rangle Iu \|_{S^{0}(J \times \mathbf{R}^{n})} \| |Iu_{h}|^{4/n} + |u_{h}|^{4/n} \|_{L_{t}^{\infty} L_{x}^{n/2}(J \times \mathbf{R}^{n})}$$

$$\lesssim \frac{1}{N^{\frac{4}{n}-}} \| \langle \nabla \rangle Iu \|_{S^{0}(J \times \mathbf{R}^{n})}^{1 + 4/n}.$$\vspace{5mm}

\noindent For $n = 3, 4$,

$$\| [|Iu|^{4/n} - |u|^{4/n}] (Iu) \|_{L_{t,x}^{2}(J \times \mathbf{R}^{n})} \lesssim \| Iu \|_{L_{t}^{\infty} L_{x}^{\frac{2n}{n - 2}}(J \times \mathbf{R}^{n})} \| |Iu|^{4/n} - |u|^{4/n} \|_{L_{t}^{2} L_{x}^{n}(J \times \mathbf{R}^{n})}$$

\noindent When $n = 4$, $$\| |Iu| - |u| \|_{L_{t}^{2} L_{x}^{4}(J \times \mathbf{R}^{n})} \lesssim \| Iu_{h} \|_{L_{t}^{2} L_{x}^{4}(J \times \mathbf{R}^{n})} + \| u_{h} \|_{L_{t}^{2} L_{x}^{4}(J \times \mathbf{R}^{n})} \lesssim \frac{1}{N^{1-}} \| \langle \nabla \rangle Iu \|_{S^{0}(J \times \mathbf{R}^{3})}^{2}.$$

\noindent When $n = 3$, $$\| |Iu|^{4/3} - |u|^{4/3} \|_{L_{t}^{2} L_{x}^{3}(J \times \mathbf{R}^{n})}$$ $$\lesssim (\| Iu_{h} \|_{L_{t}^{2} L_{x}^{6}(J \times \mathbf{R}^{3})} + \| u_{h} \|_{L_{t}^{2} L_{x}^{6}(J \times \mathbf{R}^{3})}) (\| u_{l} \|_{L_{t}^{\infty} L_{x}^{2}(J \times \mathbf{R}^{3})}^{1/3} + \| u_{h} \|_{L_{t}^{\infty} L_{x}^{2}(J \times \mathbf{R}^{3})}^{1/3})$$ $$\lesssim \frac{1}{N^{1-}} \| \langle \nabla \rangle Iu \|_{S^{0}(J \times \mathbf{R}^{3})}^{4/3}. \Box$$\vspace{5mm}

\noindent Now we tackle the main term.

\begin{lemma}\label{l5.3}

\begin{equation}\label{5.5}
\| I(|u|^{4/n} u) - (|u|^{4/n})(Iu) \|_{L_{t,x}^{2}(J \times \mathbf{R}^{n})} \lesssim \frac{1}{N^{\sigma}} \| \langle \nabla \rangle Iu \|_{S^{0}(J \times \mathbf{R}^{n})}^{1 + 4/n}.
\end{equation}
\end{lemma}

\noindent \emph{Proof:} Let

\begin{equation}\label{5.6}
f(t, \xi) = \int_{\xi_{2} + \xi_{3} = \xi} [m(\xi_{2} + \xi_{3}) - m(\xi_{3})] \widehat{F(u)}(t,\xi_{2}) \hat{u}(t,\xi_{3}) d\xi_{2}.
\end{equation}

\noindent As usual, make a Littlewood - Paley decomposition.\vspace{5mm}

\noindent \textbf{Case 1, $N_{2}, N_{3} << N$:} In this case $m(\xi_{2} + \xi_{3}) - m(\xi_{3}) \equiv 0$.\vspace{5mm}

\noindent For the remaining cases, to simplify notation, let $p_{1} = 2$, $q_{1} = \frac{2n}{n - 4}$, $p_{2} = \infty$, $q_{2} = \frac{n}{2}$ when $n > 4$ and $p_{1} = \infty$, $q_{1} = \frac{2n}{n - 2}$, $p_{2} = 2$, $q_{2} = n$ when $n = 3, 4$.\vspace{5mm}

\noindent \textbf{Case 2, $N_{2} \gtrsim N, N_{3} << N$:} In this case $|m(N_{2} + N_{3}) - m(N_{3})| \lesssim 1$. By Theorem $\ref{t2.1}$ and the Sobolev embedding theorem,

$$\sum_{N_{2} \gtrsim N} \| P_{N_{2}} F(u) \|_{L_{t}^{p_{2}} L_{x}^{q_{2}} (J \times \mathbf{R}^{n})} \sum_{N_{3} << N} \| P_{N_{3}} u \|_{L_{t}^{p_{1}} L_{x}^{q_{1}}(J \times \mathbf{R}^{n})}$$

$$\lesssim \ln(N) \sum_{N \lesssim N_{2}} \frac{1}{N_{2}^{\sigma}} \| \langle \nabla \rangle Iu \|_{S^{0}(J \times \mathbf{R}^{n})}^{1 + 4/n}$$

$$\lesssim \frac{1}{N^{\sigma-}} \| \langle \nabla \rangle Iu \|_{S^{0}(J \times \mathbf{R}^{n})}^{1 + 4/n}.$$\vspace{5mm}

\noindent \textbf{Case 3, $N_{2} << N$, $N_{3} \gtrsim N$:} In this case use the fundamental theorem of calculus,

$$|m(N_{2} + N_{3}) - m(N_{3})| \lesssim \frac{N_{2} m(N_{3})}{N_{3}}.$$

\noindent Again by Theorem $\ref{t2.1}$, Sobolev embedding,

$$\sum_{N \lesssim N_{3}} \frac{m(N_{3})}{N_{3}} \| P_{N_{3}} u \|_{L_{t}^{p_{1}} L_{x}^{q_{1}}(J \times \mathbf{R}^{n})} \sum_{N_{2} << N} N_{2} \| P_{N_{2}} F(u) \|_{L_{t}^{p_{2}} L_{x}^{q_{2}}(J \times \mathbf{R}^{n})}$$

$$\lesssim \| \langle \nabla \rangle Iu \|_{S^{0}(J \times \mathbf{R}^{n})}^{1 + 4/n} \sum_{N_{3} \gtrsim N} \frac{1}{N_{3}} \sum_{N_{2} << N} \frac{N_{2}}{N_{2}^{\sigma}}$$

$$\lesssim \frac{1}{N^{\sigma}} \| \langle \nabla \rangle Iu \|_{S^{0}(J \times \mathbf{R}^{n})}^{1 + 4/n}.$$\vspace{5mm}

\noindent \textbf{Case 4, $N_{2}, N_{3} \gtrsim N$:} In this case,

\begin{equation}\label{5.7}
\| P_{\geq N} Iu \|_{L_{t}^{p_{1}} L_{x}^{q_{1}}(J \times \mathbf{R}^{n})}  \| P_{\geq N} F(u) \|_{L_{t}^{p_{2}} L_{x}^{q_{2}}(J \times \mathbf{R}^{n})}
\end{equation}

\begin{equation}\label{5.7}
\lesssim \| \nabla Iu \|_{S^{0}(J \times \mathbf{R}^{n})}  \| P_{\geq N} F(u) \|_{L_{t}^{p_{2}} L_{x}^{q_{2}}(J \times \mathbf{R}^{n})}
\end{equation}

\begin{equation}\label{5.8}
\lesssim \frac{1}{N^{\sigma}} \| \langle \nabla \rangle Iu \|_{S^{0}(J \times \mathbf{R}^{n})}^{1 + 4/n}.
\end{equation}

\noindent Therefore, it remains to tackle $\| I(|u|^{4/n} u) \|_{L_{t,x}^{2}(J \times \mathbf{R}^{n})}$.\vspace{5mm}

\noindent \textbf{Case 4(a), $N_{2} >> N_{3}$:} In this case $m(N_{2} + N_{3}) \sim m(N_{2})$.

$$\sum_{N_{2} \gtrsim N} m(N_{2}) \| P_{N_{2}} F(u) \|_{L_{t}^{p_{2}} L_{x}^{q_{1}}(J \times \mathbf{R}^{n})} \sum_{N \lesssim N_{3} \lesssim N_{2}} \| P_{N_{3}} u \|_{L_{t}^{p_{1}} L_{x}^{q_{1}}(J \times \mathbf{R}^{n})}$$

$$\lesssim \| \langle \nabla \rangle Iu \|_{S^{0}(J \times \mathbf{R}^{n})}^{1 + 4/n} \sum \frac{m(N_{2})}{N_{2}^{\sigma s} N^{\sigma (1 - s)}} \sum_{N \lesssim N_{3} << N_{2}} \frac{N_{3}^{1 - s}}{N^{1 - s}}$$

$$\lesssim \frac{1}{N^{\sigma}} \| \langle \nabla \rangle Iu \|_{S^{0}(J \times \mathbf{R}^{n})}^{1 + 4/n}.$$\vspace{5mm}

\noindent \textbf{Case 4(b), $N_{2} << N_{3}$:} In this case $m(N_{2} + N_{3}) \sim m(N_{3})$. Choose $g(t,x)$ such that $\| g(t, x) \|_{L_{t,x}^{2}(J \times \mathbf{R}^{n})} = 1$. Then decompose

\begin{equation}\label{5.9}
\int \hat{g}(t,\xi) \hat{f}(t,\xi) d\xi.
\end{equation}

$$\sum_{N_{1} \sim N_{3}} \| P_{N_{1}} g \|_{L_{t,x}^{2}(J \times \mathbf{R}^{n})} \| P_{N_{3}} Iu \|_{L_{t}^{p_{1}} L_{x}^{q_{1}}(J \times \mathbf{R}^{n})} \sum_{N \lesssim N_{3} << N_{2}} \| P_{N_{2}} F(u) \|_{L_{t}^{p_{2}} L_{x}^{q_{1}}(J \times \mathbf{R}^{n})}$$

$$\lesssim \| \langle \nabla \rangle Iu \|_{S^{0}(J \times \mathbf{R}^{n})}^{4/n} \sum_{N_{1} \sim N_{3}} \| P_{N_{1}} g \|_{L_{t,x}^{2}(J \times \mathbf{R}^{n})} \| P_{N_{3}} \nabla Iu \|_{S^{0}(J \times \mathbf{R}^{n})} \sum_{N \lesssim N_{2} << N_{3}} \frac{1}{N_{2}^{\sigma s} N^{\sigma(1 - s)}}$$

$$\lesssim \frac{1}{N^{\sigma}} \| \langle \nabla \rangle Iu \|_{S^{0}(J \times \mathbf{R}^{n})}^{1+ 4/n},$$

\noindent using Cauchy-Schwartz.\vspace{5mm}

\noindent \textbf{Case 4(c), $N_{2} \sim N_{3}$:} In this case use the Sobolev estimate

\begin{equation}\label{5.10}
\| I(|u|^{4/n} u) \|_{L_{t,x}^{2}(J \times \mathbf{R}^{n})} \lesssim \| \nabla I(|u|^{4/n} u) \|_{L_{t}^{2} L_{x}^{\frac{2n}{n + 2}}(J \times \mathbf{R}^{n})}.
\end{equation}

\noindent In this case $|N_{2} + N_{3}| m(N_{2} + N_{3}) \lesssim m(N_{3}) |N_{3}|$.

$$\sum_{N \lesssim N_{2} \sim N_{3}} \| P_{N_{2}} F(u) \|_{L_{t}^{\infty} L_{x}^{n/2}(J \times \mathbf{R}^{n})} \| P_{N_{3}} \nabla Iu \|_{L_{t}^{2} L_{x}^{\frac{2n}{n - 2}}(J \times \mathbf{R}^{n})}$$

$$\lesssim \| \langle \nabla \rangle Iu \|_{S^{0}(J \times \mathbf{R}^{n})}^{1 + 4/n} \sum_{N \lesssim N_{2}} \frac{1}{N_{2}^{\sigma s}} \frac{1}{N^{\sigma(1 - s)}} \lesssim \frac{1}{N^{\sigma}} \| \langle \nabla \rangle Iu \|_{S^{0}(J \times \mathbf{R}^{n})}^{1 + 4/n}.$$

\noindent This completes the proof of Lemma $\ref{l5.3}$, and consequently Theorem $\ref{t5.1}$. $\Box$

\section{Proof for $n \geq 4$}
The interaction Morawetz estimates will be stated without proof.

\begin{theorem}\label{t6.1}
Suppose u solves $(\ref{0.1})$, then

\begin{equation}\label{6.1}
\| u \|_{L_{t}^{2(n - 1)} L_{x}^{\frac{2(n - 1)}{n - 2}}(J \times \mathbf{R}^{n})} \lesssim \| u_{0} \|_{L^{2}(\mathbf{R}^{n})}^{1/2} \| u \|_{L_{t}^{\infty} \dot{H}_{x}^{1/2}(J \times \mathbf{R}^{n})}^{\frac{n - 2}{n - 1}}.
\end{equation}

\noindent In addition, suppose $J = [0, T]$,

\begin{equation}\label{6.2}
\| u \|_{L_{t}^{\frac{4(n - 1)}{n}} L_{x}^{\frac{2(n - 1)}{n - 2}}(J \times \mathbf{R}^{n})} \lesssim T^{\frac{n - 2}{4(n - 1)}}\| u_{0} \|_{L^{2}(\mathbf{R}^{n})}^{1/2} \| u \|_{L_{t}^{\infty} \dot{H}_{x}^{1/2}(J \times \mathbf{R}^{n})}^{\frac{n - 2}{n - 1}}.
\end{equation}
\end{theorem}

\noindent \emph{Proof:} See $\cite{CKSTT2}$ for $n = 3$, $\cite{TVZ}$ for $n \geq 4$.\vspace{5mm}

\noindent A local well-posedness result is also needed.

\begin{theorem}\label{t6.2}
There exists $\epsilon > 0$ such that if

\begin{equation}\label{6.3}
\| u \|_{L_{t}^{\frac{4(n - 1)}{n}} L_{x}^{\frac{2(n - 1)}{n - 2}(J \times \mathbf{R}^{n})}} < \epsilon,
\end{equation}

\noindent and $\| \nabla Iu_{0} \|_{L^{2}(\mathbf{R}^{n})} \leq 1$, then

\begin{equation}\label{6.4}
\| \langle \nabla \rangle Iu \|_{S^{0}(J \times \mathbf{R}^{n})} \lesssim 1.
\end{equation}
\end{theorem}

\noindent \emph{Proof:} Let $J = [a, b]$. A solution to $(\ref{0.1})$ satisfies Duhamel's formula,

\begin{equation}\label{6.5}
Iu(t,x) = I e^{i(t - a) \Delta} u(a) + \int_{a}^{t} e^{i(t - \tau) \Delta} I(|u(\tau)|^{4/n} u(\tau)) d\tau.
\end{equation}

\noindent Make the Strichartz estimates,

$$\| \langle \nabla \rangle Iu \|_{S^{0}(J \times \mathbf{R}^{n})} \lesssim \| \langle \nabla \rangle Iu_{0} \|_{L^{2}(\mathbf{R}^{n})} + \| \langle \nabla \rangle Iu \|_{L_{t}^{2} L_{x}^{\frac{2n}{n - 2}}(J \times \mathbf{R}^{n})} \| |u|^{4/n} \|_{L_{t}^{n - 1} L_{x}^{\frac{n(n - 1)}{2(n - 2)}}(J \times \mathbf{R}^{n})}$$

$$\lesssim \| \langle \nabla \rangle Iu_{0} \|_{L^{2}(\mathbf{R}^{n})} + \epsilon^{4/n} \| \langle \nabla \rangle Iu \|_{S^{0}(J \times \mathbf{R}^{n})}.$$

\noindent So by the continuity method, for $\epsilon > 0$ sufficiently small,

\begin{equation}\label{6.6}
\| \langle \nabla \rangle Iu \|_{S^{0}(J \times \mathbf{R}^{n})} \lesssim 1.
\end{equation}

\noindent $\Box$\vspace{5mm}

\noindent \emph{Proof of Theorem $\ref{t0.3}$:} $$\int |\nabla Iu_{0}(x)|^{2} dx \leq N^{2(1 - s)} \| u_{0} \|_{\dot{H}^{s}(\mathbf{R}^{n})}^{2}.$$ By the Sobolev embedding theorem, $H^{\frac{n}{n + 4}, 2}(\mathbf{R}^{n}) \subset L^{2 + 4/n}(\mathbf{R}^{n})$, so for $n \geq 4$, $$\int |Iu_{0}(x)|^{2 + 4/n} dx \leq C \| u_{0} \|_{H^{s}(\mathbf{R}^{n})}^{2 + 4/n}.$$

\noindent Next, fix an interval $[0, T_{0}]$. Rescaling,

\begin{equation}\label{6.8}
\| u_{0, \lambda}(x) \|_{\dot{H}^{s}(\mathbf{R}^{n})} = \lambda^{-s} \| u_{0} \|_{\dot{H}^{s}(\mathbf{R}^{n})}.
\end{equation}

\noindent Therefore, choose $\lambda = C(\| u_{0} \|_{H^{s}(\mathbf{R}^{n})}) N^{\frac{1 - s}{s}}$ such that $E(Iu_{0, \lambda}) \leq \frac{1}{2}$. This also proves $|\lambda^{2} T_{0}| \lesssim N^{\alpha}$.\vspace{5mm}

\noindent Define a set

\begin{equation}\label{6.8.1}
W  = \{ t \in [0, \lambda^{2} T_{0}] : E(Iu_{\lambda}(t)) \leq \frac{3}{4} \}.
\end{equation}

\noindent We aim to prove $W = [0, \lambda^{2} T_{0}]$ for $s > \frac{n - 2}{n}$. Now $0 \in W$ and $W$ is closed, so it suffices to show $W$ is open in $[0, \lambda^{2} T_{0}]$. Suppose $W = [0, T] \subset [0, \lambda^{2} T_{0}]$, then there exists $\delta > 0$ such that $E(Iu_{\lambda}(t)) \leq 1$ on $[0, T + \delta]$.\vspace{5mm}

\noindent Next, apply the Morawetz estimates.

\begin{equation}\label{6.9}
\| P_{\leq N} u(t) \|_{\dot{H}^{1/2}(\mathbf{R}^{n})} \leq \| P_{\leq N} u(t) \|_{\dot{H}^{1}(\mathbf{R}^{n})}^{1/2} \| u_{0} \|_{L^{2}(\mathbf{R}^{n})}^{1/2}.
\end{equation}

\begin{equation}\label{6.10}
\| P_{> N} u(t) \|_{\dot{H}^{1/2}(\mathbf{R}^{n})} \leq \| P_{> N} u(t) \|_{\dot{H}^{s}(\mathbf{R}^{n})}^{1/2s} \| u_{0} \|_{L^{2}(\mathbf{R}^{n})}^{1 - 1/2s}.
\end{equation}

\noindent So if $m_{0} = \| u_{0} \|_{L^{2}(\mathbf{R}^{n})}$, then combining Theorem $\ref{t6.1}$,properties of the I - operator, and $E(Iu_{\lambda}(t)) \leq 1$ on $[0, T + \delta]$;

\begin{equation}\label{6.11}
\| u(t) \|_{L_{t}^{2(n - 1)} L_{x}^{\frac{2(n - 1)}{n - 2}}([0, T + \delta] \times \mathbf{R}^{n})} \leq C(m_{0}).
\end{equation}

\noindent Then by $(\ref{6.2})$,

\begin{equation}\label{6.12}
\| u(t) \|_{L_{t}^{\frac{4(n - 1)}{n}} L_{x}^{\frac{2(n - 1)}{n - 2}}(J \times \mathbf{R}^{n})} \lesssim \lambda^{\frac{n - 2}{2(n - 1)}} T_{0}^{\frac{n - 2}{4(n - 1)}} C(m_{0}).
\end{equation}

\noindent Partition $[0, T + \delta]$ into $$\frac{[C \lambda^{\frac{n - 2}{2(n - 1)}} T_{0}^{\frac{n - 2}{4(n - 1)}}]^{\frac{4(n - 1)}{n}}}{\epsilon^{\frac{4(n - 1)}{n}}} \sim N^{\frac{2(n - 2)(1 - s)}{ns}} T_{0}^{\frac{(n - 2)(1 - s)}{ns}}$$ subintervals with $\| u \|_{L_{t}^{\frac{4(n - 1)}{n}} L_{x}^{\frac{2(n - 1)}{n - 2}}(J_{k} \times \mathbf{R}^{n})} \leq \epsilon$ on each subinterval. Combining this with the estimate for the energy increment,

$$|E(Iu_{\lambda}(t))| \leq \frac{1}{2} + C N^{\frac{2(n - 2)}{n} \frac{1 - s}{s}} \cdot N^{-\frac{4}{n}+} T_{0}^{\frac{(n - 2)(1 - s)}{ns}}.$$

\noindent When $s > \frac{n - 2}{n}$, choosing $N$ sufficiently large, $$E(Iu_{\lambda}(t)) \leq \frac{3}{4}.$$ Therefore, $W$ is both open and closed in $[0, \lambda^{2} T_{0}]$, and $W = [0, \lambda^{2} T_{0}]$.\vspace{5mm}

\noindent \textbf{Remark:} It suffices to choose $$N \geq (4C)^{\frac{ns}{2(ns - (n - 2))}+} (T_{0}^{\frac{(n - 2)(1 - s)}{ns}})^{\frac{ns}{2(ns - (n - 2))}+}.$$

\noindent Since $\lambda = C_{0} N^{\frac{1 - s}{s}}$ and $$\| u(t) \|_{\dot{H}^{s}(\mathbf{R}^{n})} = \lambda^{s} \| u_{\lambda}(t) \|_{\dot{H}^{s}(\mathbf{R}^{n})},$$

\begin{equation}\label{6.13}
\sup_{[0, T_{0}]} \| u(t) \|_{\dot{H}^{s}(\mathbf{R}^{n})} \leq C' T_{0}^{\frac{(n - 2)(1 - s)^{2}}{2(ns - (n - 2))}},
\end{equation}

\noindent and the proof is complete. $\Box$

\section{Almost Morawetz Estimate for $n = 3$}
\noindent Since we wish to prove global well-posedness for $s > \frac{2}{5}$, it is not enough to use $$\| u \|_{L_{t,x}^{4}(J \times \mathbf{R}^{3})} \lesssim \| u \|_{L_{t}^{\infty} \dot{H}^{1/2}(J \times \mathbf{R}^{3})}^{2} \| u_{0} \|_{L^{2}(\mathbf{R}^{3})}^{2}.$$ Instead it is necessary to use almost Morawetz estimates in $n = 3$. (See \cite{CGT}, \cite{CR}, \cite{D} for a discussion of the $n = 2$ case, \cite{DPST1} for the $n = 1$ case.)

\begin{theorem}\label{t7.1}
Suppose $u$ solves the equation

\begin{equation}\label{7.0}
i u_{t} + \Delta u = |u|^{4/3} u,
\end{equation}

\noindent Suppose also that $[0, T] = \cup_{k = 1}^{K} J_{k}$.

\begin{equation}\label{7.0.1}
\| Iu \|_{L_{t,x}^{4}([0,T] \times \mathbf{R}^{3})}^{4} \lesssim \| u_{0} \|_{L^{2}(\mathbf{R}^{3})}^{3} \| Iu \|_{L_{t}^{\infty} \dot{H}^{1}(\mathbf{R}^{3})} + \frac{1}{N^{1-}} \sum_{k = 1}^{K} \| \nabla Iu \|_{S^{0}(J_{k} \times \mathbf{R}^{3})}^{16/3}.
\end{equation}

\end{theorem}

\noindent \emph{Proof:} Suppose $v$ satisfies the partial differential equation

\begin{equation}\label{7.1}
i v_{t} + \Delta v = F.
\end{equation}

\noindent Let

\begin{equation}\label{7.2}
T_{0j} = 2 Im (\overline{v(t,z)} \partial_{j} v(t,z)),
\end{equation}

\begin{equation}\label{7.3}
L_{jk} = -\partial_{j} \partial_{k} (|v(t,z)|^{2}) + 4 Re(\partial_{j} \overline{v(t,z)} \partial_{k} v(t,z)),
\end{equation}

\begin{equation}\label{7.4}
\aligned
\partial_{t} T_{0j} + \partial_{k} L_{jk} = 2(\overline{F(t,z)} \partial_{j} v(t,z) - \overline{v(t,z)} \partial_{j} F(t,z) \\ + F(t,z) \partial_{j} \overline{v(t,z)} - v(t,z) \partial_{j} \overline{F(t,z)}),
\endaligned
\end{equation}

\noindent Let $v(t,z)$ be the solution on $\mathbf{R}^{3} \times \mathbf{R}^{3}$ given by the coordinates $(x,y) = z$, $$v(t,z) = Iu(t,x) Iu(t,y).$$ Then $v(t,z)$ solves the equation

\begin{equation}\label{7.5}
\aligned
i \partial_{t} v(t,z) + \Delta_{z} v(t,z) = I(|u(t,x)|^{4/3} u(t,x)) Iu(t,y) \\ + Iu(t,x) I(|u(t,y)|^{4/3} u(t,y)).
\endaligned
\end{equation}

\noindent Next define the Morawetz action,

\begin{equation}\label{7.6}
M_{a}^{\otimes_{2}}(t) = \int_{\mathbf{R}^{3} \times \mathbf{R}^{3}} \partial_{j} a(z) \cdot Im(v(t,z) \partial_{j} \overline{v(t,z)}) dz,
\end{equation}

\noindent with $a(z) = |x - y|$.

$$\partial_{t} M_{a}^{\otimes_{2}}(t) = \int_{\mathbf{R}^{3} \times \mathbf{R}^{3}} \partial_{j} a(z) \partial_{t} T_{0j}(t,z) dz$$

$$ = -\int [\partial_{k} L_{jk}(t,z)] \partial_{j} a(z) dz $$ $$+ 2 \int [\overline{F(t,z)} \partial_{j} v(t,z) - \overline{v(t,z)} \partial_{j} F(t,z)  + F(t,z) \partial_{j} \overline{v(t,z)} - v(t,z) \partial_{j} \overline{F(t,z)}] \partial_{j} a(z) dz.$$

\begin{equation}\label{7.7}
\partial_{t} M_{a}^{\otimes_{2}}(t) = \int_{\mathbf{R}^{3} \times \mathbf{R}^{3}} \partial_{jkk} (|v(t,z)|^{2}) \partial_{j} a(z) dz
\end{equation}

\begin{equation}\label{7.8}
- 4 \int \partial_{k} Re (\partial_{j} \overline{v(t,z)} \partial_{k} v(t,z)) \partial_{j} a(z) dz
\end{equation}

\begin{equation}\label{7.9}
\aligned
+ 2 \int [\overline{F(t,z)} \partial_{j} v(t,z) - \overline{v(t,z)} \partial_{j} F(t,z) \\ + F(t,z) \partial_{j} \overline{v(t,z)} - v(t,z) \partial_{j} \overline{F(t,z)}] \partial_{j} a(z) dz.
\endaligned
\end{equation}

\noindent Integrating by parts three times in $(\ref{7.7})$, and using the identity $\Delta \Delta |x - y| = -\delta(|x - y|)$ in $\mathbf{R}^{3}$,

\begin{equation}\label{7.10}
\aligned
\int_{0}^{T} \int_{\mathbf{R}^{3} \times \mathbf{R}^{3}} \partial_{jkk} (|v(t,z)|^{2}) \partial_{j} a(z) dz \\
= \int_{0}^{T} \int_{\mathbf{R}^{3} \times \mathbf{R}^{3}} \delta(|x - y|) |Iu(t,x)|^{2} |Iu(t,y)|^{2} dx dy dt \\
= \int_{0}^{T} \int_{\mathbf{R}^{3}} |Iu(t,x)|^{4} dx dt.
\endaligned
\end{equation}

\noindent Integrating $(\ref{7.8})$ by parts once,

\begin{equation}\label{7.11}
\aligned
- 4 \int_{0}^{T} \int_{\mathbf{R}^{3} \times \mathbf{R}^{3}} \partial_{k} Re (\partial_{j} \overline{v(t,z)} \partial_{k} v(t,z)) \partial_{j} a(z) dz \\
= 4 \int_{0}^{T} \int_{\mathbf{R}^{3} \times \mathbf{R}^{3}} Re (\partial_{j} \overline{v(t,z)} \partial_{k} v(t,z)) \partial_{jk} a(z) dz.
\endaligned
\end{equation}

$$\partial_{jk}(|x - y|) = \frac{\delta_{jk}}{|x - y|} - \frac{(x - y)_{j}(x - y)_{k}}{|x - y|^{3}}.$$

\noindent This quantity is a positive definite matrix. $$\frac{\delta_{jk}}{|x - y|} - \frac{(x - y)_{j}(x - y)_{k}}{|x - y|^{3}} v_{j} v_{k} = \frac{|v|^{2}}{|x - y|} - \frac{(v \cdot (x - y))^{2}}{|x - y|^{3}} \geq 0.$$ This proves in particular that the quantity $(\ref{7.11})$ is $\geq 0$.\vspace{5mm}

\noindent Split $F = F_{g} + F_{b}$, with

\begin{equation}\label{7.12}
\aligned
F_{g} &= |Iu(t,x)|^{4/3} Iu(t,x) Iu(t,y) + |Iu(t,y)|^{4/3} Iu(t,x) Iu(t,y). \\
F_{b} &= F - F_{g}.
\endaligned
\end{equation}

\noindent Now for $(\ref{7.9})$, without loss of generality let $j = 1, 2,3$.

$$|Iu(t,y)|^{4/3} \overline{Iu(t,y) Iu(t,x)} \partial_{j} (Iu(t,x) Iu(t,y))$$ $$+ |Iu(t,y)|^{4/3} Iu(t,y) Iu(t,x) \partial_{j} (\overline{Iu(t,x)} \overline{Iu(t,y)}) = \partial_{j} (|Iu(t,y)|^{10/3} |Iu(t,x)|^{2}).$$

\noindent This cancels with the term

$$-\overline{Iu(t,x) Iu(t,y)} \partial_{j} (|Iu(t,y)|^{4/3} Iu(t,y) Iu(t,x))$$ $$- Iu(t,x) Iu(t,y) \partial_{j} (|Iu(t,y)|^{4/3} \overline{Iu(t,y) Iu(t,x)}) = -\partial_{j}(|Iu(t,y)|^{10/3} |Iu(t,x)|^{2}).$$

\noindent On the other hand,

$$|Iu(t,x)|^{4/3} \overline{Iu(t,x) Iu(t,y)} \partial_{j} (Iu(t,x) Iu(t,y))$$

$$ + |Iu(t,x)|^{4/3} (Iu(t,x) Iu(t,y)) \partial_{j} \overline{Iu(t,x) Iu(t,y)}$$

$$- \overline{Iu(t,x) Iu(t,y)} \partial_{j}(|Iu(t,x)|^{4/3} Iu(t,x) Iu(t,y))$$

$$ - Iu(t,x) Iu(t,y) \partial_{j} (|Iu(t,x)|^{4/3} \overline{Iu(t,x) Iu(t,y)})$$

$$ = -2 |Iu(t,x)|^{2} |Iu(t,y)|^{2} \partial_{j} |Iu(t,x)|^{4/3} = \frac{-4}{5} \partial_{j}(|Iu(t,y)|^{2} |Iu(t,x)|^{10/3}). $$

\noindent Integrating by parts,

$$-\int_{0}^{T} \int \partial_{j} a(z) \partial_{j} (|Iu(t,y)|^{2} |Iu(t,x)|^{10/3}) dx dy dt$$ $$= \int_{0}^{T} \int (\partial_{jj} a(z)) |Iu(t,y)|^{2} |Iu(t,x)|^{10/3} dx dy dt \geq 0.$$

\noindent All this together proves $(\ref{7.9})$ with $F$ replaced by $F_{g}$ is $\geq 0$.\vspace{5mm}

\noindent To evaluate $(\ref{7.9})$ with $F$ replaced by $F_{b}$, there are terms of the form

\begin{equation}\label{7.13}
\aligned
\int_{J} \int \partial_{j} a(x - y) |Iu(t,y)|^{2} (\partial_{j} Iu(t,x)) \\ \times [|Iu(t,x)|^{4/3} Iu(t,x) - I(|u(t,x)|^{4/3} u(t,x))] dx dy dt,
\endaligned
\end{equation}

\noindent terms of the form

\begin{equation}\label{7.14}
\aligned
\int_{J} \int \partial_{j} a(x - y) |Iu(t,y)|^{2} Iu(t,x) \\ \times \partial_{j} [|Iu(t,x)|^{4/3} Iu(t,x) - I(|u(t,x)|^{4/3} u(t,x))] dx dy dt,
\endaligned
\end{equation}

\noindent and also terms of the form

\begin{equation}\label{7.15}
\aligned
\int_{J} \int \partial_{j} a(x - y) Iu(t,x) (\partial_{j} \overline{Iu(t,x)}) Iu(t,y) \\ \times [|Iu(t,y)|^{4/3} Iu(t,y) - I(|u(t,y)|^{4/3} u(t,y))] dx dy dt.
\endaligned
\end{equation}

\noindent To evaluate a term of the form $(\ref{7.13})$, let $u_{l} = P_{\leq N/10} u$ and $u_{l} + u_{h} = u$. Then $Iu_{l} = u_{l}$, and

$$\| I(|u_{l}|^{4/3} u_{l}) - |u_{l}|^{4/3} u_{l} \|_{L_{t}^{2} L_{x}^{6/5}(J \times \mathbf{R}^{3})} \lesssim \| P_{> N/10} |u_{l}|^{4/3} u_{l} \|_{L_{t}^{2} L_{x}^{6/5}(J \times \mathbf{R}^{3})}$$

$$\lesssim \frac{1}{N} \| \nabla (|u_{l}|^{4/3} u_{l}) \|_{L_{t}^{2} L_{x}^{6/5}} \lesssim \frac{1}{N} \| u_{l} \|_{L_{t}^{\infty} L_{x}^{2}(J \times \mathbf{R}^{3})}^{4/3} \| \nabla u_{l} \|_{L_{t}^{2} L_{x}^{6}(J \times \mathbf{R}^{3})} \lesssim \frac{1}{N} \| \nabla u \|_{S^{0}(J \times \mathbf{R}^{3})}.$$

\noindent Also,

$$||u|^{4/3} u - |u_{l}|^{4/3} u_{l}| \lesssim |u_{h}| |u_{l}|^{4/3} + |u_{h}|^{7/3},$$

\noindent so

\begin{equation}\label{7.15.1}
\| I(|u|^{4/3} u) - |I u|^{4/3} (Iu) \|_{L_{t}^{2} L_{x}^{6/5}(J \times \mathbf{R}^{3})} \lesssim \frac{1}{N} \| \langle \nabla \rangle Iu \|_{S^{0}(J \times \mathbf{R}^{3})}^{7/3}.
\end{equation}

\noindent For a term of the form $(\ref{7.14})$, integrate by parts. Then $(\ref{7.14})$ is equal to an integral of the form $(\ref{7.13})$, as well as a term of the form

\begin{equation}\label{7.16}
\aligned
\int_{0}^{T} \int \partial_{jj} a(x - y) |Iu(t,y)|^{2} Iu(t,x) \\
\times [|Iu(t,x)|^{4/3} Iu(t,x) - I(|u(t,x)|^{4/3} u(t,x))] dx dy dt.
\endaligned
\end{equation}

$$\int |Iu(t,y)|^{2} \partial_{jj} a(x - y) dy$$ is controlled by a term of the form $$\int \frac{1}{|x - y|} |Iu(t,y)|^{2} dy \lesssim \| \langle \nabla \rangle Iu(t,y) \|_{L_{t}^{\infty} L_{x}^{2}(\mathbf{R}^{3})}^{2}.$$

\noindent When $|x - y| > 1$, $|\Delta a(x - y)| \lesssim 1$, so

\begin{equation}\label{7.16.1}
\sup_{t,x} \int_{|x - y| > 1} |Iu(t,y)|^{2} \Delta a(x - y) dy \lesssim \| Iu \|_{L_{t}^{\infty} L_{x}^{2}(J \times \mathbf{R}^{3})}^{2}.
\end{equation}

\noindent On $|x - y| \leq 1$,

\begin{equation}\label{7.16.2}
\sup_{t, x} \int \frac{|Iu(t,y)|^{2}}{|x - y|} dy \lesssim \| Iu \|_{L_{t}^{\infty} L_{x}^{6}(J \times \mathbf{R}^{3})}^{2} \lesssim \| \langle \nabla \rangle Iu \|_{L_{t}^{\infty} L_{x}^{2}(J \times \mathbf{R}^{3})}^{2}.
\end{equation}

\noindent Therefore, $$(\ref{7.16}) \lesssim \| |Iu(t,x)|^{4/3} Iu(t,x) - I(|u(t,x)|^{4/3} u(t,x)) \|_{L_{t}^{2} L_{x}^{6/5}(J \times \mathbf{R}^{3})}$$ $$\times \| Iu(t,x) \|_{L_{t}^{2} L_{x}^{6}(J \times \mathbf{R}^{3})} \| \langle \nabla \rangle Iu \|_{S^{0}(J \times \mathbf{R}^{3})}^{2}$$

$$\lesssim \frac{1}{N} \| \nabla Iu \|_{S^{0}(J_{k} \times \mathbf{R}^{3})}^{16/3}.$$

\noindent Finally,

$$\int_{J} \int \partial_{j} a(x,y) Iu(t,x) (\overline{\partial_{j} Iu(t,x)}) Iu(t,y) [|Iu(t,y)|^{4/3} (Iu(t,y)) - I(|u(t,y)|^{4/3} u(t,y))]$$

$$\lesssim \| \langle \nabla \rangle Iu \|_{L_{t}^{\infty} L_{x}^{2}(J \times \mathbf{R}^{3})}^{2} \| Iu \|_{L_{t}^{2} L_{x}^{6}(J \times \mathbf{R}^{3})} \| I(|u|^{4/3} u) - |Iu|^{4/3} (Iu) \|_{L_{t}^{2} L_{x}^{6/5}(J \times \mathbf{R}^{3})}$$

$$\lesssim \frac{1}{N} \| \langle \nabla \rangle Iu \|_{S^{0}(J \times \mathbf{R}^{3})}^{16/3}.$$

\noindent Putting this all together,

\begin{equation}\label{7.17}
\int_{0}^{T} \int_{\mathbf{R}^{3}} |Iu(t,x)|^{4} dx dt \lesssim |M_{a}^{\otimes_{2}}(T) - M_{a}^{\otimes_{2}}(0)| + \sum_{J_{k}} \frac{\| \langle \nabla \rangle Iu \|_{S^{0}(J_{k} \times \mathbf{R}^{3})}}{N^{1-}}.
\end{equation}

\begin{equation}\label{7.18}
\aligned
 M_{a}^{\otimes_{2}}(t)| = |\int_{\mathbf{R}^{3} \times \mathbf{R}^{3}} Iu(t,x) Iu(t,y) \partial_{j}(Iu(t,x) Iu(t,y)) dx dy| \\ \lesssim \| \nabla Iu \|_{L_{t}^{\infty} L_{x}^{2}} \| Iu \|_{L_{t}^{\infty} L_{x}^{2}}^{3},
 \endaligned
 \end{equation}

\noindent and Theorem $\ref{t7.1}$ is proved. $\Box$

\section{Proof for $n = 3$}
\noindent First prove a local well-posedness result.

\begin{theorem}\label{t8.0}
There exists $\epsilon > 0$ such that if $\| Iu \|_{L_{t}^{8/3} L_{x}^{4}(J_{k} \times \mathbf{R}^{3})} \leq \epsilon$, $\| \nabla Iu_{0} \|_{L^{2}(\mathbf{R}^{3})} \leq 1$, then $(\ref{0.1})$ has a local solution with

\begin{equation}\label{8.0}
\| \langle \nabla \rangle Iu \|_{S^{0}(J \times \mathbf{R}^{n})} \lesssim 1.
\end{equation}
\end{theorem}

\noindent \emph{Proof:} $$\| \langle \nabla \rangle Iu \|_{S^{0}(J \times \mathbf{R}^{3})} \lesssim \| \langle \nabla \rangle Iu_{0} \|_{L^{2}(\mathbf{R}^{n})} + \| \langle \nabla \rangle Iu \|_{L_{t}^{2} L_{x}^{6}(J \times \mathbf{R}^{3})} \| u \|_{L_{t}^{8/3} L_{x}^{4}(J \times \mathbf{R}^{n})}^{4/3}.$$

$$\| u \|_{L_{t}^{8/3} L_{x}^{4}(J \times \mathbf{R}^{3})} \leq \| P_{\leq N} u \|_{L_{t}^{8/3} L_{x}^{4}(J \times \mathbf{R}^{3})} + \| P_{> N} u \|_{L_{t}^{8/3} L_{x}^{4}(J \times \mathbf{R}^{3})}$$

$$ \leq \| Iu \|_{L_{t}^{8/3} L_{x}^{4}(J \times \mathbf{R}^{3})} + \frac{1}{N} \| \langle \nabla \rangle Iu \|_{L_{t}^{8/3} L_{x}^{4}(J \times \mathbf{R}^{3})}.$$

\begin{equation}\label{8.0.1}
\| \langle \nabla \rangle Iu \|_{S^{0}(J \times \mathbf{R}^{3})} \lesssim \| \langle \nabla \rangle Iu_{0} \|_{L^{2}(\mathbf{R}^{3})} + \| \langle \nabla \rangle Iu \|_{S^{0}(J \times \mathbf{R}^{3})}(\epsilon + \frac{1}{N} \| \langle \nabla \rangle Iu \|_{S^{0}(J \times \mathbf{R}^{3})})^{4/3}.
\end{equation}

\noindent Therefore, by the continuity method, $\| \langle \nabla \rangle Iu \|_{S^{0}(J \times \mathbf{R}^{3})} \lesssim 1$. $\Box$\vspace{5mm}

\noindent \emph{Proof of Theorem $\ref{t0.4}$:} In this case, $$\int |\nabla Iu_{0}(x)|^{2} dx \leq N^{2(1 - s)} \| u_{0}(x) \|_{\dot{H}^{s}(\mathbf{R}^{3})}^{2}.$$ $$\int |Iu(x)|^{10/3} dx \leq \| Iu(x) \|_{\dot{H}^{3/5}(\mathbf{R}^{3})}^{10/3} \leq N^{2 - \frac{10s}{3}} \| u(x) \|_{\dot{H}^{s}(\mathbf{R}^{n})}^{10/3}.$$

\noindent Once again, choose $\lambda = C(\| u_{0} \|_{H^{s}(\mathbf{R}^{n})}) N^{\frac{1 - s}{s}}$ such that $E(Iu_{0, \lambda}(x)) \leq \frac{1}{2}$. Define the set,

\begin{equation}\label{8.1}
W = \{ t : E(Iu_{\lambda}(t)) \leq \frac{3}{4} \} \subset [0, \lambda^{2} T_{0}].
\end{equation}

\noindent $0 \in W$, $W$ is closed. Suppose $W = [0, T]$, there exists $\delta > 0$ with $E(Iu_{\lambda}(t)) \leq 1$ on $[0, T + \delta]$.

\begin{lemma}\label{l8.2}
If $E(Iu_{\lambda}(t)) \leq 1$ on $[0, T + \delta]$,

\begin{equation}\label{8.2}
\| Iu_{\lambda}(t) \|_{L_{t,x}^{4}([0, T + \delta])}^{4} \leq \frac{3C}{2}m_{0}^{3}.
\end{equation}
\end{lemma}

\noindent \emph{Proof:} Let $$\tau = \sup \{ \tilde{T} : \| Iu_{\lambda} \|_{L_{t,x}^{4}([0, \tilde{T}] \times \mathbf{R}^{3})}^{4} \leq \frac{3C m_{0}^{3}}{2} \}.$$ If $\tau < T + \delta$, then

\begin{equation}\label{8.3}
\| Iu_{\lambda}(t) \|_{L_{t,x}^{4}([0, \tau])}^{4} \leq \frac{3C}{2} m_{0}^{3},
\end{equation}

\noindent and there exists some $\delta' > 0$ such that

\begin{equation}\label{8.4}
\| Iu_{\lambda}(t) \|_{L_{t,x}^{4}([0, \tau + \delta'])}^{4} \leq 2Cm_{0}^{3}.
\end{equation}

\noindent Then $[0, \tau + \delta']$ can be partitioned into $$\lesssim (\tau + \delta')^{1/3} \| Iu_{\lambda} \|_{L_{t,x}^{4}([0, \tau + \delta] \times \mathbf{R}^{3})}^{4} \leq (2C)^{8/3} \lambda^{2/3} T_{0}^{1/3} m_{0}^{8}$$ subintervals, and on each subinterval $J_{k}$, $\| Iu_{\lambda} \|_{L_{t}^{8/3} L_{x}^{4}(J_{k} \times \mathbf{R}^{3})} \leq \epsilon$,

\begin{equation}\label{8.5}
\| \langle \nabla \rangle Iu \|_{S^{0}(J_{k} \times \mathbf{R}^{3})} \leq C'.
\end{equation}

\noindent Applying the almost Morawetz estimate of the previous section,

\begin{equation}\label{8.6}
\| Iu_{\lambda} \|_{L_{t,x}^{4}([0, \tau + \delta'])}^{4} \leq C m_{0}^{3} + C' \frac{\lambda^{2/3} T_{0}^{1/3} (2C)^{8/3} m_{0}^{8}}{N^{1-}},
\end{equation}

\noindent so for $s > \frac{2}{5}$, take $N \gtrsim T_{0}^{\frac{s}{5s - 2}+}$,

\begin{equation}\label{8.7}
\| Iu_{\lambda} \|_{L_{t,x}^{4}([0, \tau + \delta'])}^{4} \leq \frac{3C}{2} m_{0}^{3}.
\end{equation}

\noindent Therefore, $\tau = T + \delta$. $\Box$\vspace{5mm}

\noindent Returning to the proof of the theorem,

\begin{equation}\label{8.7.1}
\| Iu_{\lambda} \|_{L_{t}^{8/3} L_{x}^{4}([0, T + \delta] \times \mathbf{R}^{3})} \lesssim (\lambda^{2} T_{0})^{1/8} \| Iu_{\lambda} \|_{L_{t,x}^{4}([0, T + \delta] \times \mathbf{R}^{3})}.
\end{equation}

\noindent Therefore, $[0, T + \delta]$ can be partitioned into $\lesssim \lambda^{2/3} T_{0}^{1/3}$ subintervals $J_{k}$ with $\| Iu_{\lambda} \|_{L_{t}^{8/3} L_{x}^{4}(J_{k} \times \mathbf{R}^{3})} \leq \epsilon$.

$$E(Iu_{\lambda}(t)) \leq \frac{1}{2} + \frac{C \lambda^{2/3} T_{0}^{1/3}}{N^{1-}}.$$

\noindent Again, choosing some $N \gtrsim T_{0}^{\frac{s}{5s - 2}+}$ with a possibly bigger constant,

$$E(Iu_{\lambda}(t)) \leq \frac{3}{4}.$$

\noindent This implies $W = [0, \lambda^{2} T_{0}]$. It suffices to take $N \gtrsim T_{0}^{\frac{s}{5s - 2}+}$, so

\begin{equation}\label{8.8}
\| u(t) \|_{{H}^{s}(\mathbf{R}^{3})} \lesssim T_{0}^{\frac{(1 - s)}{5s - 2}+}.
\end{equation}

\noindent This concludes the proof of Theorem $\ref{t0.4}$. $\Box$

\newpage

\nocite*
\bibliographystyle{plain}
\bibliography{energy}

\end{document}